\documentclass{amsart}
\usepackage{amsfonts}
\usepackage{amscd}
\usepackage{amssymb}
\usepackage{hyperref}
\newcommand{\cstar}{\bb{C}^\times}
\newcommand{\cf}{cf.\ }
\newcommand{\diag}[1]{\mbox{diag}\left\{#1\right\}}
\newcommand{\eg}{e.g.}
\newcommand{\lc}{loc.\ cit.}
\newcommand{\eqr}[1]{\mbox{(\ref{eq:#1})}}
\newcommand{\ie}{i.e.\ }
\newcommand{\Hom}{\mbox{Hom}}
\newcommand{\mult}{\mbox{\emph{Mult}}}
\newcommand{\e}[1]{\textbf{e}\left(\textstyle#1\right)}
\newcommand{\mc}[1]{\mathcal{#1}}
\newcommand{\up}{\upsilon}
\newcommand{\eq}[1]{\addtocounter{equation}{#1}(\theequation)\addtocounter{equation}{-#1}}
\newcommand{\G}{\Gamma}
\newcommand{\Teq}{T=\left(\begin{array}{cr}1&1\\0&1\end{array}\right)}
\newcommand{\Seq}{S=\left(\begin{array}{cr}0&-1\\1&0\end{array}\right)}
\newcommand{\abcd}{\begin{pmatrix}a&b\\c&d\end{pmatrix}}
\newcommand{\mat}[4]{\left(\begin{array}{rr}#1&#2\\#3&#4\end{array}\right)}
\newcommand{\twovec}[2]{\begin{pmatrix}#1\\#2\end{pmatrix}}
\newcommand{\cvec}[3]{\begin{pmatrix}#1\\#2\\#3\end{pmatrix}}
\newcommand{\bb}[1]{\mathbb{#1}}
\newcommand{\gln}[1]{GL_{#1}(\bb{C})}
\newcommand{\m}{\mc{M}}

\newtheorem{thm}{Theorem}[section]
\newtheorem{lem}[thm]{Lemma}
\newtheorem{defn}[thm]{Definition}
\newtheorem{cor}[thm]{Corollary}
\newtheorem{prop}[thm]{Proposition}

\numberwithin{equation}{section}

\begin{document}

\subjclass[2000]{11F03, 11F99}
\title[Vector-valued modular forms of dimension less than six]{Irreducible vector-valued modular forms of dimension less than six}
\author{Christopher Marks}
\address{Department of Mathematical and Statistical Sciences, University of Alberta}
\email{chris.marks@ualberta.ca}

\maketitle
\begin{abstract}
An algebraic classification is given for spaces of holomorphic vector-valued modular forms of arbitrary real weight and multiplier system,
associated to irreducible, $T$-unitarizable representations of the full modular group, of dimension less than six. For
representations of dimension less than four, it is shown that the associated space of vector-valued modular forms is a
cyclic module over a certain skew polynomial ring of differential operators. For dimensions four and five, a complete list of possible
Hilbert-Poincar\'e series is given, using the fact that the space of vector-valued modular forms is a free module over the ring
of classical modular forms for the full modular group. A mild restriction is then placed on the class of representation considered
in these dimensions, and this again yields an explicit determination of the associated Hilbert-Poincar\'e series.
\end{abstract}
\tableofcontents

\section{Introduction}
The general theory of vector-valued modular forms is now well-established in the literature, largely due to the efforts of Knopp/Mason
\cite{KM1,KM2,KM3,M1,M2} and Bantay/Gannon \cite{Bantay1,BG3,BG1,BG2}. The present paper builds upon the foundation of the Knopp/Mason theory, and
is in some sense a direct generalization, to higher dimension, of the main results of \cite{M2}. Specifically, we generalize
Theorem 5.5, \lc, which gives an algebraic classification of $\bb{Z}$-graded spaces of holomorphic vector-valued modular forms associated to two-dimensional irreducible, $T$-unitarizable representations of $\G=SL_2(\bb{Z})$. This is accomplished in \cite{M2} by establishing
that each such space is a cyclic module over a certain skew polynomial ring $\mc{R}$ of differential operators \eqr{R}, and a free module of rank two over the ring $\m=\bb{C}[E_4,E_6]$ of holomorphic, integral weight modular forms for $\G$. This classification (and in particular the computation of
the minimal weight associated to the space) is made possible by exploiting the theory of \emph{modular differential equations}, as introduced in \cite{M1}, together with Theorem 3.1 of \cite{M2}, which classifies indecomposable, $T$-semi-simple representations $\rho:\G\rightarrow\gln{2}$,
according to the eigenvalues of $\rho(T)$; here $\Teq$.

More recently, it has been shown \cite[Thm 1]{MM}, \cite[Thm 3.13]{KM3} that this free $\m$-module structure is realized in arbitrary dimension,
for an even broader class of representation than was treated in \cite{M2}. Consequently, giving an algebraic classification of spaces of vector-valued modular forms in arbitrary dimension is equivalent to determining the weights of the free generators for the $\m$-module structure of
the given spaces, including the all-important minimal weight. Furthermore, Theorem 3 of \cite{MM} gives an equivalence between the cyclicity of an $\mc{R}$-module of vector-valued modular forms, and the existence of a certain monic modular differential equation. This is significant because cyclic $\mc{R}$-modules exhibit the simplest $\m$-module structure possible and, even more importantly, the minimal weight can be determined explicitly in these cases.

Because of these advances, it is quite natural to try and generalize the techniques used in \cite{M2} to arbitrary dimension (in fact,
the main results of this paper were established \emph{before} \cite{MM} was written (\cf \cite{Marks1}) and formed the initial evidence
which led to Theorems 1 and 3 in \cite{MM}). What one requires for this task is a higher dimension analogue of Theorem 3.1 of \cite{M2},
\ie one needs to classify indecomposable, $T$-unitarizable representations of $\G$ in arbitrary dimension. Unfortunately, very little is
known about the representation theory of $\G$, even in the irreducible setting. A notable exception is \cite{TW}, which classifies
irreducible representations of the braid group $B_3$, of dimension less than six. As is well-known, $PSL_2(\bb{Z})=\G/\{\pm I\}$ is
isomorphic to the quotient of $B_3$ by its (infinite cyclic) center, and the Main Theorem of \cite{TW}, when translated into the
modular setting (Theorem \ref{thm:tw} below), serves as the desired generalization of \cite[Thm 3.1]{M2}. This result is,
to our knowledge, the strongest such generalization that exists in the literature and, furthermore (as can be seen
in the Appendix below), dimension five is in any event a natural boundary for the applicability of the techniques used here and
in \cite{M2}. It should be noted that \cite{TW} does not address the classification of indecomposable representations of $\G$,
and this creates an obstruction to using the techniques of \cite{M2} (see comments following Theorem \ref{thm:tw} below). For this
reason, we define in Section \ref{sec:dim4} a slightly restricted class of irreducible representations of $\G$, for which the theory
of modular differential equations may be applied with impunity; \cf Definition \ref{def:Tdet} below.

The layout of the paper is quite simple. In Section \ref{sec:recall}, we define the relevant terms and review the theory of
vector-valued modular forms and modular differential equations. We then proceed in subsequent Sections with the algebraic
classifications, on a dimension-by-dimension basis. Section \ref{sec:dim2} contains a quick and easy proof of \cite[Thm 5.5]{M2},
made possible by the results of \cite{MM}, and requiring no knowledge of the representation theory of $\G$, nor of the theory of
modular differential equations. The dimension three setting is handled in a completely analogous way, and this is the content of Section \ref{sec:dim3}; in particular (Theorem \ref{thm:dim3} below), we prove there that every irreducible, $T$-unitarizable $\rho:\G\rightarrow\gln{3}$
yields a space of holomorphic vector-valued modular forms which is cyclic as $\mc{R}$-module. In Sections \ref{sec:dim4} and
\ref{sec:dim5}, we first use the Free Module Theorem \cite[Thm 1]{MM} to determine the possible Hilbert-Poincar\'e series
for $\m$-modules of vector-valued modular forms of dimension four and five, respectively, and then by restricting slightly to
the \emph{$T$-determined} representations (\cf Definition \ref{def:Tdet} below), we are able to give an explicit classification
in these dimensions as well. Finally, we include an Appendix, containing what we find to be an interesting example from the theory
of modular differential equations; among other things, this example gives another indication of the impossibility of generalizing the results of \cite{TW} to dimension greater than five.

\section{Preliminaries}\label{sec:recall}
Let $\rho:\G\rightarrow\gln{d}$ denote a $d$-dimensional representation of $\G=SL_2(\bb{Z})$, $k\in\bb{R}$ an arbitrary real number, and $\up$ a multiplier system in weight $k$ (see Subsection \ref{subsec:ms} below). A function
\begin{equation}\label{eq:F}
F(z)=\cvec{f_1(z)}{\vdots}{f_d(z)}
\end{equation}
from the complex upper half-plane $\bb{H}$ to $\bb{C}^d$ is a \textbf{holomorphic vector-valued modular form} of weight $k$ (for the
 pair $(\rho,\up)$) if the following conditions are satisfied:
\begin{enumerate}
\item Each component function $f_j:\bb{H}\rightarrow\bb{C}$ is holomorphic in $\bb{H}$, and is of \emph{moderate growth at infinity},
\ie there is an integer $N\geq0$ such that $|f(x+iy)|<y^N$ holds for any fixed $x$ and $y\gg0$.\\
\item For each $\gamma\in\G$, $F|_k^\up\gamma=\rho(\gamma)F$.
\end{enumerate}
Here $|_k^\up$ denotes the standard ``slash'' action of $\G$ on the space $\mc{H}$ of holomorphic functions
$f:\bb{H}\rightarrow\bb{C}$, where for $\gamma=\abcd\in\G$, $z\in\bb{H}$ we have
\begin{equation}\label{eq:slash}
f|_k^\up\gamma(z)=\up^{-1}(\gamma)(cz+d)^{-k}f\left(\frac{az+b}{cz+d}\right).
\end{equation}
We write $\mc{H}(k,\rho,\up)$ for the $\bb{C}$-linear space of weight $k$ vector-valued modular forms for $(\rho,\up)$. If
$\rho(-I)$ is a scalar matrix, then the space $\mc{H}(\rho,\up)$ of holomorphic vector-valued
modular forms for $(\rho,\up)$ is $\bb{Z}$-graded as
\begin{eqnarray}\label{eq:zgrd}
\mc{H}(\rho,\up)=\bigoplus_{k\geq0}\mc{H}(k_0+2k,\rho,\up),
\end{eqnarray}
for some minimal weight $k_0$ which is congruent$\pmod{\bb{Z}}$ to the cusp parameter of $\up$ (see Subsection \ref{subsec:ms} below),
and satisfies the inequality $k_0\geq 1-d$ (see Corollary \ref{cor:univbound} below). Note that if $\rho(-I)$ is
\emph{not} a scalar matrix, then $\rho$ necessarily decomposes into a direct sum $\rho_+\oplus\rho_-$ of sub-representations
such that $\rho_\pm(-I)=\pm I$ (\cf comments following Lemma 2.3 in \cite{KM2}), so the assumption is merely one
of convenience.

For $U\in\gln{d}$, denote by $\rho_U$ the representation $\rho_U(\gamma)=U\rho(\gamma)U^{-1}$. As usual, we say that
$\rho$ and $\rho':\G\rightarrow\gln{d}$ are \emph{equivalent}, and write $\rho\sim\rho'$, if $\rho'=\rho_U$ for some $U$.
It is clear that in this case there is a graded isomorphism of $\bb{C}$-linear spaces
\begin{eqnarray}\label{eq:isom}
\mc{H}(\rho,\up)&\cong&\mc{H}(\rho',\up),\\
F\in\mc{H}(k,\rho,\up)&\leftrightarrow&UF\in\mc{H}(k,\rho',\up).\nonumber
\end{eqnarray}
This isomorphism allows us to focus, within a given equivalence class, on representations with particularly nice properties. Along
these lines, note that we consider in this paper only those $\rho$ which are \emph{$T$-unitarizable}, meaning that
$\rho(T)$, $\Teq$, is similar to a unitary matrix. By the above isomorphism we may, and henceforth shall, assume that
\begin{eqnarray}\label{eq:rhoT}
\rho(T)=\diag{\e{r_1},\cdots,\e{r_d}},\hspace{.3cm}0\leq r_j<1
\end{eqnarray}
(Here and throughout, we write the exponential of a real number $r$ as $\e{r}:=e^{2\pi ir}$.) Assuming this form for $\rho(T)$
ensures that the components of any $F\in\mc{H}(\rho,\up)$ have $q$-expansions familiar from the classical theory of modular
forms. In other words, slashing $F$ with the $T$ matrix and using the assumption of moderate growth shows that each component
of $F$ has the form
\begin{equation}\label{eq:qexpn}
f_j(z)=q^{\lambda_j}\sum_{n\geq0}a_j(n)q^n,
\end{equation}
where $q=e^{2\pi iz}$, for each $j$ we have
\begin{eqnarray}\label{eq:admissible}
0\leq\lambda_j\equiv r_j+\frac{m}{12}\pmod{\bb{Z}},
\end{eqnarray}
and $m$ denotes the cusp parameter of $\up$ (Subsection \ref{subsec:ms} below). We define an \textbf{admissible set}
for $(\rho,\up)$ to be any real numbers $\{\lambda_1,\cdots,\lambda_d\}$ which satisfy \eq{0}. Thus the set of
leading exponents of the components of any nonzero $F\in\mc{H}(\rho,\up)$ will, by definition, form an admissible set for
$(\rho,\up)$ (but not conversely, \ie we are not claiming that \emph{every} admissible set appears as the set of leading
exponents of some $F$). Among all admissible sets, there is a unique one which additionally satisfies $\lambda_j<1$ for each
$j$; we will refer to this as the \textbf{minimal admissible set} for $(\rho,\up)$. If $\{\lambda_1,\cdots,\lambda_d\}$ is
the minimal admissible set for $(\rho,\up)$, then every nonzero vector $F\in\mc{H}(\rho,\up)$ will have the form
\begin{eqnarray}\label{eq:arbF}
\cvec{q^{\lambda_1+n_1}\sum_{n\geq0}a_1(n)q^n}{\vdots}{q^{\lambda_d+n_d}\sum_{n\geq0}a_d(n)q^n},
\end{eqnarray}
where $a_j(0)\neq0$ for each $j$, and the $n_j$ are nonnegative integers.

We write $\m=\bigoplus_{k\geq0}\m_{2k}$ for the graded ring of integral weight, holomorphic modular forms for $\G$,
\ie for each $k\geq0$ we have $\m_{2k}=\mc{H}(2k,\textbf{1},\textbf{1})$, where $\textbf{1}$ denotes the trivial
one-dimensional representation/multiplier system, which satisfies $\textbf{1}(\gamma)=1$ for each $\gamma\in\G$. As is
well-known, $\m=\bb{C}[E_4,E_6]$ is also a graded polynomial algebra, where for each even integer $k\geq2$ we write
$$E_k(q)=1-\frac{2k}{B_k}\sum_{k\geq1}\sigma_{k-1}(n)q^n$$
for the normalized Eisenstein series in weight $k$; here $B_k$ denotes the $k^{th}$ Bernoulli number and
$\sigma_k(n)=\sum_{0<d|n}d^k$. Componentwise multiplication makes $\mc{H}(\rho,\up)$ a graded left $\m$-module, and it is
clear that the isomorphism \eqr{isom} is one of graded $\m$-modules as well as vector spaces. Regarding this structure,
one of the most important results we use in this paper (\cite[Thm 1]{MM}) is
\begin{thm}\label{thm:free}
Suppose $\rho:\G\rightarrow\gln{d}$ is $T$-unitarizable, such that $\rho(-I)$ is a scalar matrix, and let $\up$ be any
multiplier system for $\G$. Then $\mc{H}(\rho,\up)$ is a free $\m$-module of rank $d$.\qed
\end{thm}
Theorem \ref{thm:free} implies that the data needed to describe the $\m$-module structure of the graded space
\eqr{zgrd} boils down to the minimal weight $k_0$, together with $d$ distinguished nonnegative integers $k_1,\cdots,k_d$,
which give the weights $k_0+2k_1,\cdots,k_0+2k_d$ of the free generators for $\mc{H}(\rho,\up)$. Because the generators of $\m$
as graded polynomial algebra are of weights four and six, each space $\mc{H}(\rho,\up)$ has a \emph{Hilbert-Poincar\'e series}
(\cf \cite{Be}) of the form
\begin{eqnarray}\label{eq:hpseries}
\Psi(\rho,\up)(t)&=&\sum_{k\geq0}\dim\mc{H}(k_0+2k,\rho,\up)t^{k_0+2k}\nonumber\\
\nonumber\\
&=&\frac{t^{k_0}(t^{2k_1}+\cdots+t^{2k_d})}{(1-t^4)(1-t^6)}.
\end{eqnarray}
Ideally, one would like to be able to determine explicitly the Hilbert-Poincar\'e series of $\mc{H}(\rho,\up)$ for
a representation $\rho:\G\rightarrow\gln{d}$ of arbitrary dimension $d$, in terms of some invariants attached to the equivalence
class of $\rho$. It seems that the crucial step in solving this problem is the determination of the minimal weight $k_0$. For
example, it follows from the bounds developed in the proof of \cite[Thm 1]{MM} that if $\rho$ is irreducible of dimension $d$
and $\mc{H}(\rho,\up)$ has minimal weight $k_0$, then the weights of the $d$ free generators of $\mc{H}(\rho,\up)$ as $\m$-module
must lie in the interval $[k_0,k_0+2(d-1)]$; in particular, there are only a finite number of possible Hilbert-Poincar\'e series
which could describe $\mc{H}(\rho,\up)$.

Unfortunately, it is not known how to determine the minimal weight of a graded space \eqr{zgrd} for arbitrary $\rho$. In fact,
it is not even known whether the minimal weight of $\mc{H}(\rho,\up)$ has a universal upper bound as a function of $\dim\rho$, although
it seems likely that this is true (and in particular, when $\rho$ is unitary this universal bound does exist, a fact which is
implicit in the proof of the Main Theorem of \cite{KM1}). Note, however, that in situations where one is able to exploit the existence
of a vector-valued modular form arising from the solution space of a monic modular differential equation (Subsection \ref{sec:mlde} below),
and in particular in the case that $\mc{H}(\rho,\up)$ is a cyclic $\mc{R}$-module (\cf\ Theorem \ref{thm:cyclic} below), the
minimal weight \emph{can} be determined explicitly; this provides a strong motivation for elucidating the general theory of such equations.

\subsection{Multiplier systems for $\G$}\label{subsec:ms}
See \cite[Ch 3]{R} for a discussion of multiplier systems of arbitrary real weight. Note that (unlike \cite{R}) we do
\emph{not} assume that our multiplier systems are defined on $PSL_2(\bb{Z})$, thus we obtain twelve multiplier systems
for each weight, instead of the six described in \lc

Let $\bb{S}^1=\{z\in\bb{C}\ |\ |z|=1\}$ denote the unit circle. A \emph{multiplier system} in weight $k\in\bb{R}$
is a function $\up:\G\rightarrow\bb{S}^1$ which makes the map \eqr{slash} a right action of $\G$ on $\mc{H}$. The ratio
of any two multiplier systems of weight $k$ is a homomorphism, and in fact the group $\Hom(\G,\cstar)$ acts transitively
on the set $\mult(k)$ of multiplier systems in weight $k$. As is well-known, the commutator quotient of $\G$ is cyclic of
order 12, thus $\Hom(\G,\cstar)=\langle\chi\rangle$ is cyclic of order 12 as well, with generator $\chi$ satisfying
\begin{equation}\label{eq:chi}
\chi(T)=\e{\frac{1}{12}},\hspace{.3cm}\chi(S)=\e{-\frac{1}{4}},
\end{equation}
where $\Seq$. So we have, for example,
$$\mult(k)=\{\up_k\chi^N\ |\ 0\leq N\leq11\},$$
where $\up_k$ is the multiplier system which makes $\eta^{2k}$ a modular form of weight $k$; here
\begin{equation}\label{eq:eta}
\eta(q)=q^\frac{1}{24}\prod_{n\geq1}(1-q^n)
\end{equation}
denotes Dedekind's eta function. We have $\mult(k)=\mult(l)$ if and only if $k\equiv l\pmod{\bb{Z}}$, and in particular,
if $k\in\bb{Z}$ then $\up_k=\chi^k$ is a character of $\G$, so that $\mult(k)=\Hom(\G,\cstar)$ in this case.
For a given multiplier system $\up$, we define the \emph{cusp parameter} of $\up$ to be the unique real number $0\leq m<12$
such that $\up(T)=\e{\frac{m}{12}}$; note that this differs from the definition given in \cite{R} by a factor of 12.

We mention here that in Theorem \ref{thm:free} above and in the remainder of this Section, we state results in terms of
arbitrary real weight and multiplier system, whereas the reference given usually will contain a statement and proof in the integral weight,
trivial multiplier system setting. In all cases, an analysis of the proof shows that the appearance of a nontrivial
multiplier system is either inconsequential (as in the statement and proof of Theorem \ref{thm:free}), or that one may make
some trivial modifications to the original proof in order to obtain what is being claimed; most of these proofs are written down
explicitly in the present author's doctoral dissertation \cite{Marks1}, thus we will say nothing further regarding these modifications.

\subsection{Modular differential equations and the modular Wronskian}\label{sec:mlde}
Recall (\eg\ \cite[Ch 10]{L}) the \emph{modular derivative} in weight $k\in\bb{R}$,
$$D_k=\frac{1}{2\pi i}\frac{d}{dz}-\frac{k}{12}E_2=q\frac{d}{dq}-\frac{k}{12}E_2.$$
$D_k$ acts (componentwise) as a weight two operator on spaces of vector-valued modular forms, so that
$$F\in\mc{H}(k,\rho,\up)\mapsto D_kF\in\mc{H}(k+2,\rho,\up)$$
for any $(\rho,\up)$. This defines a weight two operator $D:\mc{H}(\rho,\up)\rightarrow\mc{H}(\rho,\up)$, which acts as $D_k$ on
$\mc{H}(k,\rho,\up)$, and is graded with respect to the $\m$-module structure of $\mc{H}(\rho,\up)$, \ie
$$(f,F)\in\m_k\times\mc{H}(l,\rho,\up)\mapsto D(fF)=FD_kf+fD_lF\in\mc{H}(k+l+2,\rho,\up).$$
For each $n\geq1$, we write $D_k^n$ for the composition
$$D_k^n=D_{k+2(n-1)}\circ\cdots\circ D_{k+2}\circ D_k.$$
An $d^{th}$ order \textbf{monic modular differential equation (MMDE)} in weight $k\in\bb{R}$ is an ordinary differential
equation in the disk $|q|<1$, of the form
\begin{equation}\label{eq:mlde}
L[f]=D_k^df+M_4D_k^{d-2}f+\cdots+M_{2(d-1)}D_kf+M_{2d}f=0,
\end{equation}
with $M_j\in\mc{M}_j$ for each $j$. When rewritten in terms of $\frac{d}{dq}$, one obtains from \eqr{mlde} an ODE
\begin{equation}\label{eq:rewrite}
q^df^{(d)}(q)+q^{d-1}g_{d-1}(q)f^{(d-1)}(q)+\cdots+g_0(q)f=0,
\end{equation}
where $f^{(n)}=\frac{d^nf}{dq^n}$ and $g_j$ is holomorphic in $|q|<1$ for each $j$ (in fact, each $g_j\in\bb{C}[E_2,E_4,E_6]$, the ring
of \emph{quasi-modular forms} for $\G$). Thus an MMDE has, at worst, $q=0$ as regular singular point, and no other singularities.
The theory (\cf \cite{H}) of such equations, due to Fuchs and Frobenius, tells us that if the \emph{indicial roots} of \eq{0} are
nonnegative real numbers $\lambda_1,\cdots,\lambda_d$, pairwise incongruent$\pmod{\bb{Z}}$, then the $d$-dimensional solution space $V$ of
\eqr{mlde} has a basis, which in this context is called a \emph{fundamental system of solutions} of \eqr{mlde}, consisting of functions of the form \eqr{qexpn}. It is clear that $V$ defines a subspace of $\mc{H}$, consisting of moderate growth functions. Furthermore, it is proven in \cite[Thm 4.1]{M1} that $V$, when viewed in this way, is invariant under the $|_k^\up$ action of $\G$ on $\mc{H}$, for any multiplier $\up$ in weight $k$. Thus MMDEs provide a rich source of vector-valued modular forms, as we record in the
\begin{thm}[Mason]\label{thm:mldevec}
Suppose that the MMDE \eqr{mlde} has real, nonnegative indicial roots $\lambda_1,\cdots,\lambda_d$, which are
pairwise incongruent$\pmod{\bb{Z}}$. Then there is a vector
$$F(z)=\cvec{q^{\lambda_1}+\sum_{n\geq1}a_1(n)q^{\lambda_1+n}}{\vdots}{q^{\lambda_d}+\sum_{n\geq1}a_d(n)q^{\lambda_d+n}},$$
whose components form a basis of the solution space $V$ of \eqr{mlde}, with the following property: given any multiplier system
$\up$ in weight $k$, there is a representation $\rho:\G\rightarrow\gln{d}$, arising from the $|_k^\up$ action of $\G$ on $V$,
such that $F\in\mc{H}(k,\rho,\up)$.\qed
\end{thm}
Note that if $m$ denotes the cusp parameter of $\up$, then any such $\rho$ will satisfy \eqr{rhoT}, where for each $j$
the relation \eqr{admissible} holds; in fact $\up(T)\rho(T)$ is the monodromy matrix for \eqr{mlde} at the regular singular point $q=0$,
relative to the ordered basis $\{f_1,\cdots,f_d\}$ of $V$.

We define an $n^{th}$ order \textbf{Eisenstein operator} (of weight $k\in\bb{R}$) to be an expression of the form
\begin{equation}\label{eq:eop}
L=D_k^n+\alpha_4E_4D_k^{n-2}+\cdots+\alpha_{2n}E_{2n},
\end{equation}
where $\alpha_j\in\bb{C}$ for each $j$. We have the
\begin{lem}\label{lem:eop}
Let $n\geq1$, and suppose $L[f]=0$ is an MMDE with $L$ the Eisenstein operator \eqr{eop}. Then the weight $k$ and
the $\alpha_j$ are uniquely determined by the indicial roots of the equation.
\end{lem}

\begin{proof}
First of all, note that a simple inductive argument shows that the operator $D_k^n$ can be written in the form
\begin{equation}\label{eq:dkn}
D_k^n=q^n\frac{d^n}{dq^n}+q^{n-1}f_{n,n-1}(q)\frac{d^{n-1}}{dq^{n-1}}+\cdots+f_{n,0}(q),
\end{equation}
where the $f_{n,j}$ are holomorphic in $|q|<1$ and, furthermore,
\begin{equation}\label{eq:wtk}
f_{n,n-1}(0)=\frac{n(5(n-1)-k)}{12}.
\end{equation}
If we rewrite the given MMDE in the form \eqr{rewrite} (replacing $d$ with $n$) then we have, in the notation of \eqr{dkn},
\begin{equation}\label{eq:gnj}
g_{n-j}=\left\{\begin{array}{ll}f_{n,n-1}&\ j=1\\ \\
f_{n,n-2}+\alpha_4E_4&\ j=2\\ \\
f_{n,n-j}+\alpha_{2j}E_{2j}+\sum_{i=2}^{j-1}\alpha_{2i}E_{2i}f_{n-i,n-j}
&\ j\geq3.\end{array}\right.
\end{equation}
Let $r_1,\cdots,r_n$ denote the indicial roots of the MMDE. The corresponding indicial polynomial factors as
\begin{equation}\label{eq:indprod}
\prod_{j=1}^n(r-r_j)=\sum_{j=0}^n(-1)^je_jr^{n-j},
\end{equation}
where $e_j$ denotes the $j^{th}$ elementary symmetric polynomial in $r_1,\cdots,r_n$. On the other hand, if for each $i\in\{0,1,\cdots,n-1\}$ we
define integers $a_j^i$ such that
$$\prod_{j=0}^i(r-j)=\sum_{j=1}^{i+1}a_j^ir^j,$$
then we may write the indicial equation as
\begin{equation}\label{eq:indrewrite}
r^n+A_{n-1}r^{n-1}+\cdots+A_1r+A_0=0,
\end{equation}
where we set $A_0=g_0(0)$, and for $1\leq j\leq n-1$ we define
\begin{equation}
A_{n-j}=a_{n-j}^{n-1}+\sum_{i=2}^{j+1}a_{n-j}^{n-i}g_{n-i+1}(0).
\end{equation}
Equating coefficients in \eqr{indprod} and \eqr{indrewrite}, we obtain for each $j$ the identity
\begin{equation}\label{eq:indfinal}
(-1)^je_j=A_{n-j}.
\end{equation}
Taking $j=1$ in \eqr{indfinal} yields
$$-(r_1+\cdots+r_n)=a_{n-1}^{n-1}+a_{n-1}^{n-2}g_{n-1}(0),$$
and combining this with \eqr{gnj} and \eqr{wtk} (and noting that $a_i^{i-1}=1$ for any $i\geq0$) produces
the formula
$$k=\frac{12}{n}\left(a_{n-1}^{n-1}+\sum_{j=1}^nr_j\right)+5(n-1),$$
thus the weight $k$ of the MMDE is determined uniquely by the $r_j$ (note that this also follows from Theorem
\ref{thm:MMDEwt} below). For $j=2$, \eqr{indfinal} and \eqr{gnj} (recall also that $E_4(0)=1$ by definition) yield
$$e_2=a_{n-2}^{n-1}+a_{n-2}^{n-2}f_{n,n-1}(0)+f_{n,n-2}(0)+\alpha_4,$$
so by \eqr{wtk} we have $\alpha_4$ as a function of $k$ and the indicial roots; since we have just shown that
$k$ itself is a function of the $r_j$, we see that $\alpha_4$ is as well. For arbitrary $j\geq3$, \eqr{gnj} and \eqr{indfinal} show that $\alpha_{2j}$ is a function of the indicial roots and $k,\alpha_4,\cdots,\alpha_{2(j-1)}$. If we assume inductively that $k$ and $\alpha_{2i}$, $2\leq i\leq j-1$ are determined uniquely by the indicial roots of the MMDE, then we find that $\alpha_{2j}$  is as well.
\end{proof}

\begin{cor}\label{cor:ummde} Let $n\leq5$. For each set $\{\lambda_1,\cdots,\lambda_n\}$
of complex numbers, there is a unique $n^{th}$ order MMDE with indicial roots $\lambda_1,\cdots,\lambda_n$.
\end{cor}

\begin{proof}
This follows directly from the previous Lemma and the fact that $\m_{2j}$ is spanned by $E_{2j}$ for $j=2,3,4,5$,
so that every MMDE of order less than 6 is of the form $L[f]=0$, with $L$ an Eisenstein operator.
\end{proof}

As is well-known, the solution space of the first-order MMDE $D_kf=0$ is spanned by $\eta^{2k}$, and this immediately implies
\begin{lem}
Assume that $\dim\rho\geq2$, and that $F\in\mc{H}(\rho,\up)$ has linearly independent components. Then $DF\neq0$.\qed
\end{lem}
Thus we have the useful
\begin{cor}\label{cor:Dinj}
If $\rho$ is irreducible, $\dim\rho\geq2$, then for each multiplier system $\up$, $D$ is an injective operator on $\mc{H}(\rho,\up)$.\qed
\end{cor}
We will also make use of
\begin{lem}\label{lem:dfind} Suppose $F=(f_1,\cdots,f_d)^t\in\mc{H}(k,\rho,\up)$ has
linearly independent components. Then for each $n\leq d$, the set $\{F,D_kF,\cdots,D_k^{n-1}F\}$ is independent
over $\m$; in particular, this set spans a rank $n$ free submodule
$$\bigoplus_{j=0}^{n-1}\m D_k^jF$$
of $\mc{H}(\rho,\up)$.
\end{lem}

\begin{proof}
Suppose there is a relation
\begin{equation}
M_{n-1}D_k^{n-1}F+M_{n-2}D_k^{n-2}F+\cdots+M_1D_kF+M_0F=0
\end{equation}
with $M_j\in\m$ for each $j$. Rewriting \eq{0} in terms of $\frac{d}{dq}$ yields an ordinary differential equation
$L[f]=0$ of order at most $n-1$, for which each of the $d$ linearly independent components of $F$ is a solution. By the
Fuchsian theory of ODEs in the complex domain, this is impossible unless $L$ is identically 0,
and one sees easily that this forces $M_j=0$ for each $j$.
\end{proof}

The \textbf{modular Wronskian}, defined in \cite[Sec 3]{M1}, plays a key role in the techniques we use in the current paper.
We gather here various results (Thms 3.7, 4.3 \lc) in the following
\begin{thm}[Mason]\label{thm:MMDEwt}
Assume $F\in\mc{H}(k,\rho,\up)$ is of the form \eqr{arbF}, and has linearly independent components; set
$\lambda=\sum\lambda_j$, $n=\sum n_j$. Then the modular Wronskian of $F$ has the form
$$W(F)=\eta^{24(\lambda+n)}g\in\mc{H}(d(d+k-1),\det\rho,\up^d),$$
for some nonzero modular form $g\in\m_{d(d+k-1)-12\lambda}$ which is not a cusp form. In particular,
the weight $k$ of $F$ satisfies the inequality
\begin{equation}\label{eq:wtbound}
k\geq\frac{12(\lambda+n)}{d}+1-d,
\end{equation}
and equality holds in \eq{0} if, and only if, the components of $F$ span the solution space of an MMDE
\eqr{mlde} in weight $k$.\qed
\end{thm}

Since the exponents $\lambda_j+n_j$ in \eqr{arbF} are nonnegative, we obtain from \eq{0} a universal lower bound on the
minimal weight $k_0$ in \eqr{zgrd}:

\begin{cor}\label{cor:univbound}
Assume $\rho$ is irreducible with $\rho(T)$ given by \eqr{rhoT}, let $\up$ be an arbitrary multiplier system, and let
$\{\lambda_1,\cdots,\lambda_d\}$ be the minimal admissible set for $(\rho,\up)$, with $\lambda=\sum\lambda_j$. Then the
minimal weight $k_0$ in \eqr{zgrd} satisfies the inequality
\begin{equation}\label{eq:univbound}
k_0\geq\frac{12\lambda}{d}+1-d.
\end{equation}\qed
\end{cor}

One also obtains from the modular Wronskian the important observation that the representations ``of MMDE type'' -- \ie those
representations arising from the slash action of $\G$ on the solution space of an MMDE -- are always indecomposable:

\begin{lem}\label{lem:indecMMDE}
Suppose $\rho:\G\rightarrow\gln{d}$ is $T$-unitarizable, $\up\in\mult(k)$, and $\mc{H}(\rho,\up)$
contains a vector
\begin{equation}\label{eq:MMDEvec}
F=\cvec{f_1}{\vdots}{f_d}=\cvec{q^{\lambda_1}+\sum_{n\geq1}a_1(n)q^n}{\vdots}{q^{\lambda_d}+\sum_{n\geq1}a_d(n)q^n}
\end{equation}
whose components form a fundamental system of solutions of an MMDE in weight $k$.  Then $\rho$ is indecomposable.
\end{lem}

\begin{proof}
Suppose $\rho$ decomposes into a direct sum $\rho=\rho_1\oplus\rho_2$. We may assume, up to equivalence of representation,
that the $|_k^\up$-invariant subspaces corresponding to $\rho_1$ and $\rho_2$ are spanned by $\{f_1,\cdots,f_{d_1}\}$, $\{f_{\scriptscriptstyle{d_1+1}},\cdots,f_d\}$ respectively, for some $1\leq d_1\leq d$. Set $\dim\rho_2=d_2=d-d_1$, and
$\Lambda_1=\lambda_1+\cdots+\lambda_{d_1}$, $\Lambda_2=\lambda_{d_1+1}+\cdots+\lambda_d$. By Theorem \ref{thm:MMDEwt}, we have
\begin{equation}
d(k+d-1)=12(\Lambda_1+\Lambda_2).
\end{equation}
On the other hand, if we define $F_1=(f_1,\cdots,f_{d_1})^t$, $F_2=(f_{d_1+1},\cdots,f_d)^t$, then $F_1\in\mc{H}(k,\rho_1,\up)$,
$F_2\in\mc{H}(k,\rho_2,\up)$, and Theorem \ref{thm:MMDEwt} yields the inequalities
$$d_1(k+d_1-1)\geq12\Lambda_1,$$
$$d_2(k+d_2-1)\geq12\Lambda_2.$$
Adding these inequalities and using $\eq{0}$ yields the
inequality
$$2d_1d_2\leq0,$$
so that $d_1=d$, $d_2=0$.
\end{proof}

Finally, we recall (\cite[Sec 2]{M2}) the \emph{skew polynomial ring}
\begin{equation}\label{eq:R}
\mc{R}=\{f_0+f_1d+\cdots+f_nd^n\ |\ f_j\in\m,n\geq0\}
\end{equation}
of differential operators, which combines the actions of $\m$ and the modular derivative on $\mc{H}(\rho,\up)$. Addition
is defined in $\mc{R}$ as though it were the polynomial ring $\m[d]$, and multiplication is performed via the identity
$$df=fd+Df,$$
where $D$ denotes the modular derivative. Each space $\mc{H}(\rho,\up)$ of vector-valued modular forms is a $\bb{Z}$-graded
left $\mc{R}$-module in the obvious way, and again we point out that the isomorphism \eqr{isom} is one of graded $\mc{R}$-modules
as well as vector spaces. Regarding this structure, we record here another key result from \cite{MM}, which
will be used frequently in subsequent Sections:

\begin{thm}\label{thm:cyclic}
Suppose $\rho:\G\rightarrow\gln{d}$ satisfies \eqr{rhoT}, $\rho(-I)$ is a scalar matrix, and let $\up$ be any multiplier system for $\G$. Let $\{\lambda_1,\cdots,\lambda_d\}$ denote the minimal admissible set for $(\rho,\up)$,
put $\lambda=\sum_{j=1}^d\lambda_j$, and write $\mc{H}(\rho,\up)$ as in \eqr{zgrd}. Then the following hold:
\begin{enumerate}
\item If $\mc{H}(\rho,\up)=\mc{R}F_0$ is cyclic as $\mc{R}$-module, then $F_0$ has the form \eqr{MMDEvec},
and the $f_j$ form a fundamental system of solutions of a $d^{th}$-order MMDE of weight $k_0=\frac{12\lambda}{d}+1-d$. The indicial roots of the MMDE are $\lambda_1,\cdots,\lambda_d$, and they are \emph{distinct}.\\
\item Conversely, suppose that the $\lambda_j$ are distinct. Then there is a $d^{th}$-order MMDE in weight
$k_0=\frac{12\lambda}{d}+1-d$, such that
$$\mc{H}(\rho',\up)=\bigoplus_{k\geq0}\mc{H}(k_0+2k,\rho',\up)=\mc{R}F$$
is cyclic as $\mc{R}$-module; here $F$ is as in \eqr{MMDEvec}, the $f_j$ span the solution space $V$ of the MMDE, and $\rho'$ denotes the representation of $\G$ arising from the $|_{k_0}^\up$-action of $\G$ on $V$, relative to the ordered basis $\{f_1,\cdots,f_d\}$. Consequently, $\rho'(T)=\rho(T)$, and $\rho'$ is indecomposable by Lemma \ref{lem:indecMMDE}, so in particular $\rho'(-I)$ is a scalar matrix.
\end{enumerate}\qed
\end{thm}

This completes the necessary review of the basic theory of vector-valued modular forms and modular differential equations.
We now proceed to the classification of spaces $\mc{H}(\rho,\up)$ for irreducible $\rho$ of dimension $d\leq5$. As a warm-up,
we derive here the well-known results from the classical (\ie one-dimensional) setting, using the vector-valued
methods:

Fix an integer $0\leq N\leq11$ and a character $\chi^N:\G\rightarrow\bb{C}^\ast$, with $\chi$ as in \eqr{chi}, and let
$\up$ be a multiplier system for $\G$, with cusp parameter $m$. Then $\chi^N(T)=\e{\frac{N}{12}}$, and the minimal admissible
set for $(\chi^N,\up)$ is $\{\lambda_1\}$, where $\lambda_1$ satisfies the congruence \eqr{admissible}, with $r_1=\frac{N}{12}$.
Viewing $\lambda_1$ as an indicial root, we obtain by Corollary \ref{cor:ummde} a unique first order MMDE $D_{k_0}f=0$, where
by Theorem \ref{thm:MMDEwt} we have $k_0=\frac{\lambda_1}{12}$. The solution space $V$ of the MMDE is spanned by a function
$$f_1(z)=q^{\lambda_1}+\sum_{n\geq1}a(n)q^{\lambda_1+n}\in\mc{H}(k_0,\rho,\up),$$
where $\rho:\G\rightarrow\bb{C}^\ast$ is the representation afforded us by Theorem \ref{thm:mldevec}, which arises from the
$|_{k_0}^\up$ action of $\G$ on $V$, relative to the basis $\langle f_1\rangle$ of $V$. In fact $f_1=\eta^{2k_0}$, as is well-known,
and from Theorem \ref{thm:cyclic} we obtain a cyclic $\mc{R}$-module
$$\mc{H}(\rho,\up)=\bigoplus_{k\geq0}\mc{H}(k_0+2k,\rho,\up)=\mc{R}\eta^{2k_0},$$
which in the one-dimensional setting is equivalent to saying that $\mc{H}(\rho,\up)$ is a free $\m$-module of rank one, with
generator $\eta^{2k_0}$; this is the content of Theorem \ref{thm:free}, as it pertains to the present setting. Note that
by definition we have $\eta^{2k_0}|_{k_0}^\up T=\rho(T)\eta^{2k_0}$, so from \eqr{admissible} we conclude that
$\rho(T)=\e{\frac{N}{12}}$. Recalling that a character of $\G$ is completely determined by its value at the matrix $T$, this
shows that $\rho=\chi^N$, so we have classified our space $\mc{H}(\chi^N,\up)$.

It will be seen in what follows that this same method may, to some extent, be utilized in any dimension less than six.

\section{Dimension two}\label{sec:dim2}
This has been worked out in the trivial multiplier system case in \cite{M2}. Here we extend the results to arbitrary real weight,
and provide a streamlined (indeed, nearly trivial) proof, made possible by Theorems \ref{thm:free} and \ref{thm:cyclic}.
\begin{thm}\label{thm:2d}
Let $\rho:\G\rightarrow\gln{2}$ be irreducible with $\rho(T)$ as in \eqr{rhoT}, fix a multiplier system $\up$, and let
$\{\lambda_1,\lambda_2\}$ be the minimal admissible set of $(\rho,\up)$. Then
$$\mc{H}(\rho,\up)=\bigoplus_{k\geq0}\mc{H}(k_0+2k,\rho,\up)=\mc{R}F_0$$
is cyclic as $\mc{R}$-module, with $k_0=6(\lambda_1+\lambda_2)-1$, and the components of $F_0$ form a fundamental system of
solutions of a second order MMDE in weight $k_0$.
\end{thm}

\begin{proof}
Write $\mc{H}(\rho,\up)$ as in \eqr{zgrd}. It is clear that the number of weight $k_0$ generators
of $\mc{H}(\rho,\up)$ (as $\m$-module) is exactly $\dim\mc{H}(k_0,\rho,\up)$. Similarly, since $\m_2=\{0\}$, the number of
weight $k_0+2$ generators is $\dim\mc{H}(k_0+2,\rho,\up)$. But if we fix any nonzero $F_0\in\mc{H}(k_0,\rho,\up)$, then by
Corollary \ref{cor:Dinj} we know that $DF_0\in\mc{H}(k_0+2,\rho,\up)$ is nonzero, so by Theorem \ref{thm:free},
we conclude that $\mc{H}(k_0,\rho,\up)=\langle F_0\rangle$ and $\mc{H}(k_0+2,\rho,\up)=\langle DF_0\rangle$ are 1-dimensional,
and $\mc{H}(\rho,\up)=\m F_0\oplus\m DF_0$ as $\m$-module. In particular, $\mc{H}(\rho,\up)=\mc{R}F_0$ is cyclic as $\mc{R}$-module,
and part 1 of Theorem \ref{thm:cyclic} finishes the proof.
\end{proof}

\begin{cor} The Hilbert-Poincar\'e series of $\mc{H}(\rho,\up)$ is
$$\Psi(\rho,\up)(t)=\frac{t^{k_0}(1+t^2)}{(1-t^4)(1-t^6)},$$
thus for each $k\geq0$ we have (using the well-known dimension formula for $\m_k$)
$$\dim\mc{H}(k_0+2k,\rho,\up)=\left[\frac{k}{3}\right]+1.$$
\qed\end{cor}
We now give an alternate proof of Theorem \ref{thm:2d}, along the same lines as the method used to classify spaces for
one-dimensional representations at the end of the last Section. This is roughly the method used in the original proof
found in \cite{M2}, and its successful application relies on the following Theorem, \cite[Thm 3.1]{M2}:
\begin{thm}[Mason]\label{thm:mason2d} Suppose $\rho:\G\rightarrow\gln{2}$ is indecomposable, with $\rho(T)=\diag{x_1,x_2}$ for some $x_j\in\bb{C}$.
Then the following are equivalent:
\begin{enumerate}
\item $\rho$ is irreducible.\\
\item The ratio $x_1/x_2$ is \emph{not} a primitive sixth root of 1.\\
\item The eigenvalues $\{x_1,x_2\}$ of $\rho(T)$ define a unique equivalence class of 2-dimensional indecomposable representations
of $\G$.
\end{enumerate}\qed
\end{thm}
Assume once again the hypotheses of Theorem \ref{thm:2d}. By Corollary \ref{cor:ummde}, there is a unique second order MMDE
$$D^2_{k_0}f+\alpha_4E_4f=0$$
whose set of indicial roots is exactly $\{\lambda_1,\lambda_2\}$, the minimal admissible set of $(\rho,\up)$. By part two
of Theorem \ref{thm:cyclic}, we have $k_0=6(\lambda_1+\lambda_2)-1$, and there is a cyclic $\mc{R}$-module
$$\mc{H}(\rho',\up)=\mc{R}F_0=\bigoplus_{k\geq0}\mc{H}(k_0+2k,\rho',\up),$$
where $\rho':\G\rightarrow\gln{2}$ is a representation arising from the $|_{k_0}^\up$ action of $\G$ on the solution space $V$
of the MMDE, and the generator
$$F_0=\twovec{q^{\lambda_1}+\cdots}{q^{\lambda_2}+\cdots}\in\mc{H}(k_0,\rho',\up)$$
has components which span $V$. We have $\rho'(T)=\rho(T)$ and Theorem \ref{thm:cyclic} (really Lemma \ref{lem:indecMMDE}) says
that $\rho'$ is indecomposable. Applying Theorem \ref{thm:mason2d}, we see that $\rho$ and $\rho'$ are equivalent irreducible
representations, so the isomorphism \eqr{isom} establishes Theorem \ref{thm:2d}.\qed

\section{Dimension three}\label{sec:dim3}
This is completely analogous to dimension two, and we will prove quite easily
\begin{thm}\label{thm:dim3}
Let $\rho:\G\rightarrow\gln{3}$ be an irreducible representation with $\rho(T)$ as in \eqr{rhoT}, fix a multiplier system $\up$,
and write $\{\lambda_1,\lambda_2,\lambda_3\}$ for the minimal admissible set of $(\rho,\up)$. Then
$$\mc{H}(\rho,\up)=\bigoplus_{k\geq0}\mc{H}(k_0+2k,\rho,\up)=\mc{R}F_0$$
is cyclic as $\mc{R}$-module, with $k_0=4(\lambda_1+\lambda_2+\lambda_3)-2$, and the components of $F_0$ form a fundamental system of
solutions of a third order MMDE in weight $k_0$.
\end{thm}

\begin{proof}
Write $\mc{H}(\rho,\up)$ as the graded sum \eqr{zgrd}. It follows immediately from Corollary \ref{cor:Dinj} that $\dim\mc{H}(k_0,\rho,\up)=1$,
since otherwise there would be two generators $F,G$ of weight $k_0$ and two $DF,DG$ of weight $k_0+2$, in violation of
Theorem \ref{thm:free}. Fix any nonzero $F_0$ of minimal weight, and write $\mc{H}(k_0,\rho,\up)=\langle F_0\rangle$. Again
by Corollary \ref{cor:Dinj}, we may take $DF_0$ as a second generator of $\mc{H}(\rho,\up)$. Suppose there is another generator
of weight $k_0+2$, say $G$. Then by Theorem \ref{thm:free} we have $$\mc{H}(\rho,\up)=\m F_0\oplus\m DF_0\oplus\m G.$$
But then $D^2F_0\in\mc{H}(k_0+4,\rho,\up)$ must satisfy a relation $D^2F_0=M_4F_0$, with $M_4\in\m_4$. This is impossible by
Lemma \ref{lem:dfind}, so we must have $\mc{H}(k_0+2,\rho,\up)=\langle DF_0\rangle$, and $D^2F_0$ can be taken as the third
generator of $\mc{H}(\rho,\up)$. Therefore $\mc{H}(\rho,\up)=\mc{R}F_0$ is cyclic as $\mc{R}$-module, and the rest of the Theorem
follows from part 1 of Theorem \ref{thm:cyclic}.
\end{proof}

\begin{cor} The Hilbert-Poincar\'e series of $\mc{H}(\rho,\up)$ is
$$\Psi(\rho,\up)(t)=\frac{t^{k_0}(1+t^2+t^4)}{(1-t^4)(1-t^6)},$$
so that
$$\dim\mc{H}(k_0+2k,\rho,\up)=\left[\frac{k}{2}\right]+1$$
for all $k\geq0$.
\qed\end{cor}

\section{Dimension four}\label{sec:dim4}
We first record a technical result which we will need in this Section:
\begin{lem}\label{lem:div}
Let $\rho:\G\rightarrow\gln{4}$ be irreducible with $\rho(T)$ as in \eqr{rhoT}, and set $r=\sum r_j$.
Then $3r\in\bb{Z}$.
\end{lem}

\begin{proof}
Because $\rho$ is irreducible, we know that $\rho(S^2)=\rho(-I)=\pm I_4$, so the eigenvalues of $\rho(S)$
are $\pm1$, $\pm i$ respectively. Note that in either case both eigenvalues occur, since $\rho$ is irreducible and,
as is well-known, $S$ and $T$ generate $\G$. Define $R=TS^{-1}=\mat{-1}{1}{-1}{0}$. Then $R^3=I$ and $R,S$ generate $\G$
as well. If $\rho(S)$ has a three-dimensional eigenspace $U$, then the nonzero subspace $U\cap\rho(R)U\cap\rho(R^2)U$
is invariant under both $\rho(R)$ and $\rho(S)$, and this violates the irreducibility of $\rho$. Therefore the eigenvalues
of $\rho(S)$ are either $\{1,1,-1,-1\}$ or $\{i,i,-i,-i\}$, and either way we have $\det\rho(S)=1$. This implies
$$\e{3r}=\det\rho(T^3)=\det\rho(RS)^3=1,$$
so $3r$ is an integer.
\end{proof}

Continuing with the assumptions of the Lemma, fix a multiplier system $\up$, and write the minimal admissible set of $(\rho,\up)$ as $\{\lambda_1,\lambda_2,\lambda_3,\lambda_4\}$. From the relations \eqr{admissible} and Lemma \ref{lem:div},
we find that $3\lambda\equiv m\pmod{\bb{Z}}$, so the minimal weight $k_0$ in \eqr{zgrd} must be of the form
$3\lambda+N$ for some integer $N$. Furthermore, if one takes a nonzero vector $F$ of minimal weight, then the identity
$F|_{3\lambda+N}^\up S^2=\rho(S^2)F$ implies the relation
\begin{eqnarray}\label{eq:4dident}
\rho(S^2)=\up(S^2)^{-1}(-1)^{3\lambda+N}I_4.
\end{eqnarray}
In particular, if $\mc{H}(\rho,\up)$ is cyclic as $\mc{R}$-module, then by Theorem \ref{thm:cyclic} we know $k_0=3\lambda-3$,
so \eq{0} holds exactly when $N$ is odd; this provides a necessary, though perhaps not sufficient, criterion for
$\mc{H}(\rho,\up)$ to be cyclic as $\mc{R}$-module. Regardless, it turns out that the lowest weight space for
$\mc{H}(\rho,\up)$ will be one-dimensional, as we now prove:
\begin{lem}\label{lem:4dmin}
Let $F$ be an arbitrary nonzero vector in $\mc{H}(k_0,\rho,\up)$, written as in \eqr{arbF}. Then $n_j=0$ for $j=1,2,3,4$.
\end{lem}

\begin{proof}
Suppose otherwise, so that $n_i\geq1$ for some fixed $i\in\{1,2,3,4\}$, and consider the subspace
$$V=\langle E_{10}F,E_8DF,E_6D^2F,E_4D^3F\rangle\leq\mc{H}(k_0+10,\rho,\up).$$
By Lemma \ref{lem:dfind} we have $\dim V=4$, and it is clear that any nonzero $G\in V$, written in the form \eqr{arbF},
will again satisfy $n_i\geq1$. For each $j\in\{1,2,3,4\}-\{i\}$, let $\phi_j:V\rightarrow\bb{C}$ denote the linear
functional which takes such a $G$ to $\phi_j(G)=a_j(0)$, the first Fourier coefficient of the
$j^{th}$ component of $G$. Then $\dim\ker\phi_j\geq3$ for each $j$, so that
$$\bigcap_{j\neq i}\ker\phi_j\neq\{0\}.$$
This is equivalent to saying there is a nonzero $G\in V$ which satisfies $n_j\geq1$ for $j=1,2,3,4$, when written in the form \eqr{arbF}.
In other words, recalling the weight 12 cusp form
$$\Delta(q)=\eta^{24}(q)=q\prod_{n\geq1}(1-q^n)^{24}\in\m_{12},$$
we have that $\frac{G}{\Delta}$ is a nonzero vector in $\mc{H}(k_0-2,\rho,\up)$. But this cannot be, since $k_0$ is
by definition the minimal weight for $\mc{H}(\rho,\up)$, and this contradiction finishes the proof.
\end{proof}

\begin{cor}\label{cor:k01d}
$\dim\mc{H}(k_0,\rho,\up)=1$.
\end{cor}

\begin{proof}
Suppose there are two linearly independent vectors $F,G$ in $\mc{H}(k_0,\rho,\up)$. Since $F$ and $G$ each satisfy the
conclusion of Lemma \ref{lem:4dmin}, it is clear that some linear combination of these vectors will produce a nonzero vector in
$\mc{H}(k_0,\rho,\up)$ which violates the conclusion of the Lemma.
\end{proof}

With these results in hand, we are now able to show that there is only one possible non-cyclic structure for $\mc{H}(\rho,\up)$:
\begin{lem}\label{lem:4noncyclic}
If $\mc{H}(\rho,\up)$ is \emph{not} cyclic as $\mc{R}$-module, then it has the Hilbert-Poincar\'e series
\begin{equation}\label{eq:hp4}
\Psi(\rho,\up)(t)=\frac{t^{k_0}(1+2t^2+t^4)}{(1-t^4)(1-t^6)},
\end{equation}
with corresponding dimension formula
$$\dim\mc{H}(k_0+2k,\rho,\up)=\left[\frac{2k+1}{3}\right]+1,$$
for all $k\geq0$. The minimal weight $k_0$ is congruent$\pmod{\bb{Z}}$ to the cusp parameter of $\up$, and satisfies the
inequality \eqr{univbound}.
\end{lem}

\begin{proof}
Corollary \ref{cor:k01d} implies that $\mc{H}(k_0,\rho,\up)$ contributes exactly one generator to the $\m$-module structure
of $\mc{H}(\rho,\up)$, say F, and we know that $DF\in\mc{H}(k_0+2,\rho,\up)$ can be taken as a second generator. If
$\mc{H}(\rho,\up)$ is not cyclic as $\mc{R}$-module, and if there is \emph{not} a second generator of weight $k_0+2$, then
there must be two of weight $k_0+4$. But this would imply a relation $D^3F=\alpha_1E_6F+\alpha_2DF$, in violation of Lemma
\ref{lem:dfind}. Thus the Hilbert-Poincar\'e series indicated is the correct one.
\end{proof}

Thus there are exactly two possible $\m$-module structures in the four-dimensional irreducible setting. Unfortunately, in the
most general case we are not able to say definitively which of the two structures obtains, given the input $(\rho,\up)$;
neither are we able to determine explicitly the minimal weight $k_0$. We \emph{can} say the following:

\begin{cor} Suppose that \eqr{4dident} holds for $N$ even. Then the Hilbert-Poincar\'e series of $\mc{H}(\rho,\up)$ is
given by \eqr{hp4}, and $k_0=3\lambda+N$ for some even $N\geq-2$.\end{cor}

\begin{proof} By the comments following \eqr{4dident}, it is clear that the Hilbert-Poincar\'e series is the non-cyclic one, and
the inequality $N\geq-2$ follows from \eqr{univbound}.
\end{proof}

In the case where \eqr{4dident} holds for $N$ \emph{odd}, one would like to say that $\mc{H}(\rho,\up)$ is a cyclic $\mc{R}$-module,
but again, this is not known in complete generality. Nonetheless, in the vast majority of cases, we \emph{are} able to determine
quite explicitly what occurs. To see this, we will employ the method given in the alternative proof of Theorem \ref{thm:2d} above,
together with the following results concerning the representation theory of $\G$, due to Tuba and Wenzl
(\cf \cite{TW}, Corollary in Section 2, Main Theorem 2.9, and subsequent Corollary):

\begin{thm}[Tuba/Wenzl]\label{thm:tw}
Let $\rho:\G\rightarrow\gln{d}$ be an irreducible representation of dimension $d\leq5$. Then the following hold:
\begin{enumerate}
\item The minimal and characteristic polynomials of $\rho(T)$ coincide.\\
\item If $d\neq4$, then the eigenvalues of $\rho(T)$ define a unique equivalence class of irreducible representations.\\
\item If $d=4$, there are at most two equivalence classes of irreducible representations defined by the eigenvalues of $\rho(T)$.
\end{enumerate}\qed
\end{thm}

The above Theorem provides a higher-dimensional generalization of Theorem \ref{thm:mason2d}, with the significant \emph{caveat}
that it says nothing about indecomposable representations which might be lurking about, with the same eigenvalues at $T$ as the
given irreducible representation $\rho$; note that an important consequence of Theorem \ref{thm:mason2d} is that this does
\emph{not} happen in dimension two. To our knowledge, there is currently no classfication theory for indecomposable representations
of $\G$, in any dimension, apart from Theorem \ref{thm:mason2d}. Thus it is not known to what extent this phenomenon occurs in dimensions 3,4,5, \ie how often an
indecomposable-but-not-irreducible representation $\rho'$ occurs such that $\rho'(T)=\rho(T)$ for some irreducible $\rho$;
see, however, the Appendix below for an explicit example which shows that this phenomenon definitely \emph{does} occur in every
dimension greater than five. In any event, this concept presents an obstruction to the use of the MMDE theory for the classification
of spaces of four- or five-dimensional vector-valued modular forms: one may, as in the two-dimensional setting, construct MMDEs which
produce representations $\rho'$ such that $\rho'(T)=\rho(T)$, where $\rho$ is the given irreducible representation, but Lemma
\ref{lem:indecMMDE} tells us only that $\rho'$ is indecomposable, so in the most general context Theorem \ref{thm:tw} might not
apply. In order to overcome this deficiency in our method, we must restrict to those irreducible representations which have no
``shadow'' indecomposables. To this end, we make the

\begin{defn}\label{def:Tdet}
Let $d\geq1$. An irreducible representation $\rho:\G\rightarrow\gln{d}$ is \textbf{T-determined} if the following condition holds:
\begin{tabbing}\hspace{2cm}\=If $\rho':\G\rightarrow\gln{d}$ is indecomposable and $\rho'(T)$
has the same\\
\>eigenvalues as $\rho(T)$, then $\rho'$ is irreducible.
\end{tabbing}
\end{defn}

In fact, this is a very mild restriction. For example, given an arbitrary representation $\rho$, if no proper sub-product
of the eigenvalues of $\rho(T)$ is a $12^{th}$ root of 1 then $\rho$ is $T$-determined; this follows from the fact that
$\det\rho$ is a character $\chi^N$ of $\G$, with $\chi$ as in \eqr{chi}. Also, Theorem \ref{thm:mason2d} implies that every
irreducible $\rho:\G\rightarrow\gln{2}$ with $\rho(T)$ semi-simple is $T$-determined. And for unitary representations,
of course, the notions of indecomposable and irreducible coincide. In any event,
for the $T$-determined representations we are able to give a completely explicit classification in dimension four:

\begin{thm}\label{thm:4d}
Let $\rho:\G\rightarrow\gln{4}$ be a $T$-determined representation with $\rho(T)$ as in \eqr{rhoT}, fix a multiplier system
$\up$ for $\G$, let $\{\lambda_1,\cdots,\lambda_4\}$ be the minimal admissible set for $(\rho,\up)$, and set $\lambda=\sum\lambda_j$.
Then exactly one of the following holds:

\begin{enumerate}
\item The relation \eqr{4dident} holds for $N$ odd, and $\mc{H}(\rho,\up)$ is cyclic as $\mc{R}$-module, with minimal weight
$3\lambda-3$.\\
\item The relation \eqr{4dident} holds for $N$ even, and $\mc{H}(\rho,\up)$ has the Hilbert-Poincar\'e series \eqr{hp4},
with minimal weight $3\lambda-2$.
\end{enumerate}\end{thm}

\begin{proof}
Note first that, since $\rho(T)$ is diagonal, part one of Theorem \ref{thm:tw} and the relations \eqr{admissible}
imply that the $\lambda_j$ are distinct. As in the alternate proof of Theorem \ref{thm:2d}, one views the $\lambda_j$ as
the indicial roots of a unique (by Corollary \ref{cor:ummde}) MMDE of order 4. By part two of Theorem \ref{thm:cyclic},
we obtain a cyclic $\mc{R}$-module
$$\mc{H}(\rho_0,\up)=\mc{R}F_0=\bigoplus_{k\geq0}\mc{H}(3\lambda-3+2k,\rho_0,\up),$$
for some representation $\rho_0$ which satisfies $\rho_0(T)=\rho(T)$; note that \eqr{4dident} is satisfied by $\rho_0$ for
$N$ odd. On the other hand, we may take as indicial roots the set $\{\lambda_1+1,\lambda_2,\lambda_3,\lambda_4\}$, and obtain
a fourth order MMDE in weight $3(\lambda+1)-3=3\lambda$. From Theorem \ref{thm:mldevec}, we obtain a nonzero vector
\begin{equation}\label{eq:F1}
F_1=\left(\begin{array}{l}q^{\lambda_1+1}+\cdots\\q^{\lambda_2}+\cdots\\q^{\lambda_3}+\cdots\\q^{\lambda_4}+\cdots
\end{array}\right)\in\mc{H}(3\lambda,\rho_1,\up),
\end{equation}
whose components span the solution space of this second MMDE. Note that this second representation $\rho_1$ -- arising from
the $|_{3\lambda}^\up$ action of $\G$ on the solution space of the MMDE -- satisfies \eqr{4dident} for $N$ even. The isomorphism
\eqr{isom} then shows that $\rho_0$ and $\rho_1$ are inequivalent representations, each of which is indecomposable thanks to Lemma \ref{lem:indecMMDE}. Furthermore, we have $\rho_0(T)=\rho_1(T)=\rho(T)$. Since $\rho$ is $T$-determined, this forces $\rho_0$ and
$\rho_1$ to be irreducible. By part two of Theorem \ref{thm:tw}, we find that $\rho$ is equivalent to exactly one of the $\rho_j$, thus $\mc{H}(\rho,\up)$ is isomorphic to $\mc{H}(\rho_j,\up)$ for $j=0$ or 1; this can be determined explicitly, of course, by examining the
relation \eqr{4dident}. We have already seen that if \eqr{4dident} holds for $N$ odd, then $\rho$ is equivalent to $\rho_0$, and part
one of the Theorem obtains.

On the other hand, suppose \eqr{4dident} holds for $N$ even, so that $\rho$ is equivalent to $\rho_1$. Then we already know
that $\mc{H}(\rho,\up)$ has the Hilbert-Poincar\'e series \eqr{hp4}, so to finish the proof we need only determine the minimal
weight for this module. By \eqr{isom}, it suffices to do this for $\mc{H}(\rho_1,\up)$. We know there is a nonzero vector in $\mc{H}(3\lambda,\rho_1,\up)$, so by \eqr{zgrd} and Corollary \ref{cor:univbound}, the minimal weight is either $3\lambda$ or
$3\lambda-2$. But \eqr{F1} and Lemma \ref{lem:4dmin} make it clear that $3\lambda$ \emph{cannot} be the minimal weight, so it
must be $3\lambda-2$.
\end{proof}

\section{Dimension five}\label{sec:dim5}
In this Section, we again apply Theorem \ref{thm:free}, etc., to determine the possible Hilbert-Poincar\'e series for
spaces of vector-valued modular forms associated to five-dimensional irreducible representations of $\G$, and then restrict
to the \emph{$T$-determined} setting in order to obtain the most explicit results possible via our methods.

Assume for the remainder of this Section that $\rho:\G\rightarrow\gln{5}$ is irreducible with $\rho(T)$ as in \eqr{rhoT},
fix a multiplier system $\up$ with cusp parameter $m$, and write $\{\lambda_1,\cdots,\lambda_5\}$ for the minimal
admissible set of $(\rho,\up)$. If we write each $\lambda_j$ in the form
$$\lambda_j=r_j+\frac{m}{12}+l_j,\hspace{.3cm}l_j\in\{-1,0\}$$
and set $\lambda=\sum\lambda_j$, $r=\sum r_j$, $l=\sum l_j$, then we have
\begin{eqnarray*}
\frac{12\lambda}{5}-4&=&\frac{12}{5}\left(r+\frac{5m}{12}+l\right)-4\\ \\
&=&m+\frac{12(r+l)}{5}-4,
\end{eqnarray*}
so that
$$\frac{12\lambda}{5}-4\equiv m\pmod{\bb{Z}}\ \Leftrightarrow\ 12(r+l)\equiv0\pmod{5}.$$
This shows in particular (via Theorem \ref{thm:cyclic}) that if $\mc{H}(\rho,\up)$ is cyclic as $\mc{R}$-module, then
necessarily
$$\frac{12(r+l)}{5}\in\bb{Z}.$$
Certainly this does not hold in general, and instead represents a very special case. Note that $12r\in\bb{Z}$, since
$\det\rho$ is a character of $\G$, and $(5,12)=1$, so in fact there is a unique $N\in\{0,1,2,3,4\}$ such that
\begin{eqnarray}\label{eq:5cong}
12(r+l+N)\equiv0\pmod{5}.
\end{eqnarray}
Thus the minimal weight for $\mc{H}(\rho,\up)$, whatever it turns out to be, will necessarily be of the form
\begin{equation}\label{eq:5dminwt}
\frac{12(\lambda+N)}{5}-4+n,
\end{equation}
for some $n\geq-\frac{12N}{5}$ (this follows from Corollary \ref{cor:univbound}), and $N=n=0$ exactly when
$\mc{H}(\rho,\up)$ is cyclic as $\mc{R}$-module.

We will see below that in some sense, the possible Hilbert-Poincar\'e series for $\mc{H}(\rho,\up)$ correspond to the
values of $N$ in the discussion above. The following result starts us down this path, by determining explicitly these possibilities:
\begin{thm}\label{thm:5dhp}Write $\mc{H}(\rho,\up)$ as the graded $\m$-module \eqr{zgrd}. Then there are five possibilities
for the associated Hilbert-Poincar\'e series, namely
\begin{equation}\label{eq:5dhp}
\Psi(\rho,\up)(t)=\frac{t^{k_0}P_N(t)}{(1-t^4)(1-t^6)},
\end{equation}
where $k_0\geq\frac{12\lambda}{5}-4$ and
\begin{eqnarray*}
P_0(t)&=&1+t^2+t^4+t^6+t^8,\\
P_1(t)&=&2+2t^2+t^4,\\
P_2(t)&=&1+t^2+2t^4+t^6,\\
P_3(t)&=&1+2t^2+t^4+t^6,\\
P_4(t)&=&1+2t^2+2t^4.
\end{eqnarray*}
\end{thm}

\begin{proof}
The bound claimed for the minimal weight $k_0$ is just that provided by the modular Wronskian \eqr{univbound}.

Obviously the $N=0$ case occurs exactly when $\mc{H}(\rho,\up)$ is cyclic as $\mc{R}$-module, so assume for the remainder
of the proof that this is \emph{not} the case. Then at least one of the spaces $\mc{H}(k_0+2k,\rho,\up)$, $0\leq k\leq4$, contains more
than one generator for the $\m$-module structure of $\mc{H}(\rho,\up)$. On the other hand, it is clear from Theorem
\ref{thm:free} and Corollary \ref{cor:Dinj} that $\dim\mc{H}(k_0,\rho,\up)<3$, since assuming otherwise would produce at least
six generators in the $k_0,k_0+2$ spaces alone.

Suppose for the moment that $\mc{H}(k_0,\rho,\up)$ is two-dimensional. By taking an appropriate linear combination of vectors, it is clear that we may produce an $1\leq i\leq5$ and a nonzero $F\in\mc{H}(k_0,\rho,\up)$ which, when written
in the notation \eqr{arbF}, satisfies $n_i\geq1$. Then every vector in the subspace
$$V=\langle E_4D^4F,E_6D^3F,E_8D^2F,E_{10}DF,E_{12}F\rangle\leq\mc{H}(k_0+12,\rho,\up)$$
satisfies $n_i\geq1$ as well, when written as in \eqr{arbF}. We again argue via linear functionals, as in Lemma
\ref{lem:4dmin}, and (noting that $\dim V=5$ by Lemma \ref{lem:dfind}) obtain a nonzero vector $G\in\mc{H}(k_0,\rho,\up)$
such that $\Delta G\in V$. Lemma \ref{lem:dfind} and the fact that $\Delta G\in V$ make it clear that $G$ is not a
scalar multiple of $F$, so we have $\mc{H}(k_0,\rho,\up)=\langle F,G\rangle$. By Corollary \ref{cor:Dinj}, we know that
$\mc{H}(k_0+2,\rho,\up)$ is \emph{at least} two-dimensional, as it contains the subspace $\langle DF,DG\rangle$; thus four
of the five generators predicted by Theorem \ref{thm:free} are already accounted for. Note that the fifth generator 
is found in $\mc{H}(k_0+2,\rho,\up)$ if and only if $\dim\mc{H}(k_0+2,\rho,\up)=3$, and in this case we would know from Corollary \ref{cor:Dinj} that $\dim\mc{H}(k_0+4,\rho,\up)\geq3$. But this arrangement would also imply that $\mc{H}(k_0+4,\rho,\up)=\langle E_4F,E_4G\rangle$
is two-dimensional, since $\m_2=\{0\}$ and there would be no additional generators in $\mc{H}(k_0+4,\rho,\up)$. This
contradiction shows that $\mc{H}(k_0+2,\rho,\up)=\langle DF,DG\rangle$ is in fact two-dimensional.

We claim that $D^2F$ may be taken as the fifth generator for the $\m$-module structure of $\mc{H}(\rho,\up)$. To prove this, it suffices to establish that $D^2F\not\in\langle E_4F,E_4G\rangle$,
but this is clear from the fact that $\Delta G\in V$, \ie assuming a relation
$$D^2F=\alpha_1E_4F+\alpha_2E_4G$$
for some $\alpha_1,\alpha_2\in\bb{C}$ and multiplying by $\Delta$ would produce a relation which violates the
conclusion of Lemma \ref{lem:dfind}. Consequently $\mc{H}(\rho,\up)$ has the Hilbert-Poincar\'e series \eqr{5dhp}, with $N=1$, and in all the remaining cases we should, and will now, assume that $\mc{H}(k_0,\rho,\up)$ is one-dimensional:

Suppose that $\dim\mc{H}(k_0+2,\rho,\up)=1$ as well, and fix any nonzero $F$ in $\mc{H}(k_0,\rho,\up)$. Then we infer from the
hypotheses and Corollary \ref{cor:Dinj} that
$$\mc{H}(k_0,\rho,\up)=\langle F\rangle,\ \mc{H}(k_0+2,\rho,\up)=\langle DF\rangle.$$
We may also take $D^2F\in\mc{H}(k_0+4,\rho,\up)$ as a third generator, by Lemma \ref{lem:dfind}. Note that, since $\mc{H}(\rho,\up)$
is not cyclic as $\mc{R}$-module, either $\mc{H}(k_0+4,\rho,\up)$ or $\mc{H}(k_0+6,\rho,\up)$ contains more than
one generator. In particular, if
$$\mc{H}(k_0+4,\rho,\up)=\langle E_4F,D^2F\rangle$$
is two-dimensional then there must be a vector $G$ such that
$$\mc{H}(k_0+6,\rho,\up)=\langle D^3F,E_4DF,E_6F,G\rangle$$
is four-dimensional, and $D^3F,G$ can be taken as the fourth and fifth generators of $\mc{H}(\rho,\up)$.
This would mean $\mc{H}(k_0+8,\rho,\up)=\langle E_8F,E_6DF,E_4D^2F\rangle$ is three-dimensional, yet (by Corollary \ref{cor:Dinj})
contains the four-dimensional subspace $D\mc{H}(k_0+6,\rho,\up)$, contradiction. Thus $\dim\mc{H}(k_0+4,\rho,\up)\geq3$, which
implies (again by Corollary \ref{cor:Dinj}) that $\dim\mc{H}(k_0+6,\rho,\up)\geq3$. Since the three known generators $F,DF,D^2F$
only produce the two-dimensional subspace $\langle E_6F,E_4DF\rangle\leq\mc{H}(k_0+6,\rho,\up)$, it is clear in this case
that the Hilbert-Poincar\'e series for $\mc{H}(\rho,\up)$ is \eqr{5dhp} with $N=2$.

Finally, suppose that $\dim\mc{H}(k_0,\rho,\up)=1$, $\dim\mc{H}(k_0+2,\rho,\up)=2$. By the hypotheses and Corollary \ref{cor:Dinj},
there are generators $F,G$ such that $\mc{H}(k_0,\rho,\up)=\langle F\rangle$, $\mc{H}(k_0+2,\rho,\up)=\langle DF,G\rangle$.
Lemma \ref{lem:dfind}, together with the fact that $\m_2=\{0\}$, makes it clear that $D^2F$ can be taken as a fourth generator.
If $DG$ is \emph{not} contained in the subspace $\langle E_4F,D^2F\rangle\leq\mc{H}(k_0+4,\rho,\up)$, then the five generators are $F,G,DF,D^2F,DG$ and we have
the $N=4$ case of the Theorem. Otherwise, there is a relation
\begin{eqnarray}\label{eq:free1}
DG=\alpha_1D^2F+\alpha_2E_4F
\end{eqnarray}
for some $\alpha_j\in\bb{C}$, and to show that the $N=3$ case obtains, it suffices to show that the set $\{F,G,DF,D^2F,D^3F\}$
is independent over $\m$.

Assume there is a homogeneous relation in weight $k_0+2k$, say
\begin{eqnarray}\label{eq:free2}
QG=M_{2k}F+M_{2(k-1)}DF+M_{2(k-2)}D^2F+M_{2(k-3)}D^3F,
\end{eqnarray}
where $Q\in\m_{2(k-1)}$, and $M_j\in\m_j$ for each $j$. If $Q\neq0$, then dividing by $Q$ and taking the modular derivative in
\eqr{free2} yields, after utilizing \eqr{free1}, a relation
\begin{eqnarray*}
0&=&\frac{M_{2(k-3)}}{Q}D^4F
+\left[\frac{M_{2(k-2)}}{Q}+D\left(\frac{M_{2(k-3)}}{Q}\right)\right]D^3F\\ \\
&+&\left[\frac{M_{2(k-1)}}{Q}+D\left(\frac{M_{2(k-2)}}{Q}\right)-\alpha_1\right]D^2F\\ \\
&+&\left[\frac{M_{2k}}{Q}+D\left(\frac{M_{2(k-1)}}{Q}\right)\right]DF
+\left[D\left(\frac{M_{2k}}{Q}\right)-\alpha_2E_4\right]F.
\end{eqnarray*}
Noting that Lemma \ref{lem:dfind} obviously still holds when the coefficient functions lie in the fraction field of $\m$,
we conclude that each of the coefficients in the above equation is identically zero. It is then apparent that all the
$M_j$ must be zero, so that no relation like \eqr{free2} exists with a nonzero $Q$. But again by Lemma \ref{lem:dfind}, if $Q=0$ in
\eqr{free2} then all the $M_j$ are zero as well. This concludes the proof of the Theorem.
\end{proof}

As with the analogous statements in the previous Section, the above Theorem serves as an existence result only, since we
face in dimension five the same obstruction discussed after the statement of Theorem \ref{thm:tw}: the existence of indecomposable
representations which are \emph{not} irreducible, yet have the same eigenvalues at $T$ as some irreducible representation. As in
dimension four, in these cases we cannot say definitively which Hilbert-Poincar\'e series obtains in Theorem \ref{thm:5dhp}, nor
can we determine explicitly the minimal weight. However, if we again fall back into the \emph{$T$-determined} setting
(\cf Definition \ref{def:Tdet}), everything is quite explicit, as we now show. We begin with an important
\begin{lem}\label{lem:5dtdet}Retaining the hypotheses and notations from the beginning of this Section, make the additional
assumption that $\rho$ is $T$-determined, let $N\in\{0,1,2,3,4\}$ be the unique integer such that \eqr{5cong} holds,
and set 
$$k_N=\frac{12(\lambda+N)}{5}-4.$$ 
Then the following statements hold:
\begin{enumerate}
\item For each set $\{N_1,\cdots,N_5\}$ of nonnegative integers such that $\sum N_j=N$, there is a unique fifth order
MMDE with indicial roots $\lambda_j+N_j$, of weight $k_N$, and a vector-valued modular form
\begin{equation}\label{eq:5dMMDEvec}
F_{(N_1,\cdots,N_5)}(z)=\cvec{q^{\lambda_1+N_1}+\cdots}{\vdots}{q^{\lambda_5+N_5}+\cdots}\in\mc{H}(k_N,\rho',\up)
\end{equation}
whose components span the solution space of the MMDE. The representation $\rho'$ (which depends on the $N_j$) is equivalent to
$\rho$, and we have
$$\mc{H}(\rho',\up)=\bigoplus_{k\geq n_N}\mc{H}(k_N+2k,\rho',\up),$$
where the integer $n_N$ satisfies the inequality $n_N\geq-\frac{6N}{5}$.\\ \\
\item For each integer $k\geq-\frac{6N}{5}$, we have
$$\dim\mc{H}(k_N+2k,\rho',\up)\leq\left\{\begin{array}{lll}
\left[\frac{5k}{6}\right]+N&\ &k\equiv 5\pmod{6}\\ \\
\left[\frac{5k}{6}\right]+N+1&\ &k\not\equiv 5\pmod{6}.\end{array}\right.$$\\ \\
\end{enumerate}
\end{lem}

\begin{proof}
Given any appropriate set of $N_j$, the existence and uniqueness of the MMDE follows from Corollary \ref{cor:ummde}.
Note that part one of Theorem \ref{thm:tw}, the assumption that $\rho(T)$ is diagonal, and the relations \eqr{admissible}
imply that the indicial roots $\lambda_j+N_j$ are incongruent$\pmod{\bb{Z}}$, thus by Theorem \ref{thm:mldevec} we obtain
the vector $F_{(N_1,\cdots,N_5)}$ and representation $\rho'$. Noting that $\rho'(T)=\rho(T)$, we find from Lemma
\ref{lem:indecMMDE} and the assumption that $\rho$ is $T$-determined that $\rho'$ and $\rho$ are equivalent. The fact that the
minimal weight of $\mc{H}(\rho',\up)$ is congruent$\pmod{\bb{Z}}$ to $k_N$ follows from the discussion at the beginning of
this Section, and the inequality on $n_N$ follows from the bound \eqr{univbound}. This establishes part one of the
Lemma.

As for part two, suppose that $k\geq-\frac{6N}{5}$, and assume there is a nonzero vector $F$ in $\mc{H}(k_N+2k,\rho',\up)$
of the form \eqr{arbF}. By Theorem \ref{thm:MMDEwt}, the modular Wronskian of $F$ is of the form $W(F)=\eta^{24(\lambda+n)}g$,
where $n=\sum n_j$, and the weight of the non-cusp form $g$ is
\begin{eqnarray*}
wt(g)&=&5(k_N+2k+4)-12(\lambda+n)\\ \\
&=&5\left(\frac{12(\lambda+N)}{5}-4+2k+4\right)-12(\lambda+n)\\ \\
&=&12(N-n)+10k.
\end{eqnarray*}
In particular, $g$ is nonzero, so we know that $wt(g)\geq0,\neq2$, as $\m_k=\{0\}$ if $k<0,=2$.
We clearly have
$$wt(g)\equiv2\pmod{12}\ \Leftrightarrow\ k\equiv5\pmod{6},$$
and this gives the inequalities
\begin{equation}\label{eq:wbound}
n\leq\left\{\begin{array}{lll}\left[\frac{5k}{6}\right]+N-1&\ &k\equiv 5\pmod{6}\\ \\
\left[\frac{5k}{6}\right]+N&\ &k\not\equiv 5\pmod{6}.\end{array}\right.
\end{equation}
Now one may argue using linear functionals, as in the proof of Lemma \ref{lem:4dmin}, and conclude that the
$\mc{H}(k_N+2k,\rho',\up)$ will always contain a vector of the form \eqr{arbF}, such that the inequality
$$n\geq\dim\mc{H}(k_N+2k,\rho',\up)-1$$
obtains. Combining this fact with \eqr{wbound} completes the proof of the Lemma.
\end{proof}

Using this Lemma, we may now establish
\begin{thm}\label{thm:finthm} Assume the hypotheses and conclusions of Lemma \ref{lem:5dtdet}. Then the Hilbert-Poincar\'e series of
$\mc{H}(\rho,\up)$ is of the form
\begin{equation}\label{eq:hpfinal}
\Psi(\rho,\up)(t)=\Psi_N(t)=\frac{t^{k_N+n_N}P_N(t)}{(1-t^4)(1-t^6)},
\end{equation}
with $P_N$ as in the statement of Theorem \ref{thm:5dhp}, and $n_N$ given by
$$n_N=0,0,-2,-3,-4\mbox{ for }N=0,1,2,3,4$$
respectively.
\end{thm}

\begin{proof}
Thanks to Lemma \ref{lem:5dtdet} and the isomorphism \eqr{isom}, it suffices to determine the Hilbert-Poincar\'e series for
$\mc{H}(\rho',\up)$, where $\rho'$ is any representation arising from the $|_{k_N}^\up$ action of $\G$ on the solution space
of any MMDE with indicial roots $\lambda_j+N_j$ satisfying $\sum N_j=N$. We now proceed on a case-by-case basis, depending
on $N$:

$\textbf{N=0.}$ By Theorem \ref{thm:cyclic}, $\mc{H}(\rho',\up)$ is cyclic as $\mc{R}$-module, thus the Hilbert-Poincar\'e
series is of the form \eqr{hpfinal} with $N=0$, $n_0=0$, as claimed.

$\textbf{N=1.}$ Choose the integers $N_j$ as $N_1=1$, $N_j=0$ for $2\leq j\leq5$. Let $F_1=F_{(1,0,0,0,0)}\in\mc{H}(k_1,\rho',\up)$
denote the vector \eqr{5dMMDEvec}, and note that Lemma \ref{lem:5dtdet} implies that $k_1$ is the minimal weight space
for $\mc{H}(\rho',\up)$, since $-1\equiv5\pmod{6}$. We claim that the minimal weight space $\mc{H}(k_1,\rho',\up)$ is two-dimensional.
To verify this claim, we use the linear functional argument from Lemma \ref{lem:4dmin} together with Lemma \ref{lem:dfind}, and produce a vector
\begin{equation}
\widetilde{G}=\beta_1E_4D^4F_1+\beta_2E_6D^3F_1+\beta_3E_8D^2F_1+\beta_4E_{10}DF_1+\beta_5E_{12}F_1
\end{equation}
in $\mc{H}(k_1+12,\rho,\up)$ such that $G=\frac{\widetilde{G}}{\Delta}$ is a nonzero vector in $\mc{H}(k_1,\rho',\up)$.
If there is a relation $\alpha_1F_1+\alpha_2G=0$, then multiplying by $\Delta$ and substituting with \eq{0} yields
\begin{eqnarray}
0&=&\alpha_2\beta_1E_4D^4F_1+\alpha_2\beta_2E_6D^3F_1+\alpha_2\beta_3E_8D^2F_1+\\
&\ &\alpha_2\beta_4E_{10}DF_1+[\alpha_1\Delta+\alpha_2\beta_5E_{12}]F_1.\nonumber
\end{eqnarray}
By Lemma \ref{lem:dfind}, each coefficient function in \eq{0} must be zero. In particular, we have
$\alpha_1\Delta=-\alpha_2\beta_5E_{12}$. Comparing $q$-expansions forces
$\alpha_1=0$, which means $\alpha_2G=0$, \ie $\alpha_2=0$. Therefore $F$ and $G$ are linearly independent, so
by part two of Lemma \ref{lem:5dtdet} (with $(N,k)=(1,0)$), we conclude that $\mc{H}(k_1,\rho',\up)=\langle F_1,G\rangle$
is two-dimensional. Theorem \ref{thm:5dhp} then makes it clear that the Hilbert-Poincar\'e series for $\mc{H}(\rho',\up)$
must be of the form \eqr{hpfinal} with $N=1$, $n_1=0$.

$\textbf{N=2.}$ Set $F_2=F_{(1,1,0,0,0)}\in\mc{H}(k_2,\rho',\up)$ in \eqr{5dMMDEvec}. As in the $N=1$ case just treated,
we use $F_2$ to produce a nonzero vector
$$G'\in\langle D^4F_2,E_4D^2F_2,E_6DF_2,E_8F_2\rangle\leq\mc{H}(k_2+8,\rho',\up)$$
such that $G=\frac{G'}{\Delta}$ is a nonzero vector in $\mc{H}(k_2-4,\rho',\up)$. Part one of Lemma \ref{lem:5dtdet} implies
that the minimal weight of $\mc{H}(\rho',\up)$ is in fact $k_2-4$, and part two with $(k,N)=(-2,2)$ shows that
$\mc{H}(k_2-4,\rho',\up)=\langle G\rangle$ is one-dimensional. Similarly, part two of the Lemma with $(k,N)=(-1,2)$ and
Corollary \ref{cor:Dinj} show that $\mc{H}(k_2-2,\rho',\up)=\langle DG\rangle$ is one-dimensional. Since we know from
Theorem \ref{thm:cyclic} that $\mc{H}(\rho',\up)$ is not cyclic, this information is enough to conclude, via Theorem
\ref{thm:5dhp}, that the Hilbert-Poincar\'e series of $\mc{H}(\rho,\up)$ is given by \eqr{hpfinal} with $N=2$, $n_2=-2$,
as claimed.

$\textbf{N=3.}$ We shall be more particular in this case about our choice of indicial roots, singling
out a $j_1\in\{1,\cdots,5\}$ such that $\lambda_{j_1}\neq k_3-6$; this is certainly possible, since the diagonal nature
of $\rho(T)$ and part one of Theorem \ref{thm:tw} imply that the $\lambda_j$ are distinct. Having made this selection, we
then choose distinct $j_2,j_3,j_4,j_5\in\{1,\cdots,5\}-\{j_1\}$, and set $N_{j_1}=N_{j_2}=0$, $N_{j_3}=N_{j_4}=N_{j_5}=1$.
Using these integers, we then proceed as usual, and set $F_3=F_{(N_1,\cdots,N_5)}\in\mc{H}(k_3,\rho',\up)$ in \eqr{5dMMDEvec}.
We find, using the linear functional argument, a vector
$$G\in\langle D^3F_3,E_4DF_3,E_6F_3\rangle\leq\mc{H}(k_3+6,\rho',\up)$$
such that $G_1=\frac{G}{\Delta}$ is a nonzero vector in $\mc{H}(k_3-6,\rho',\up)$. Part one of Lemma \ref{lem:5dtdet}
shows that $k_3-6$ is the minimal weight for $\mc{H}(\rho',\up)$, and part two with $(k,N)=(-3,3)$ shows that
$\mc{H}(k_3-6,\rho',\up)=\langle G_1\rangle$ is one-dimensional. A second application of this reasoning produces a
vector
$$H\in\langle D^4F_3,E_4D^2F_3,E_6DF_3,E_8F_3\rangle\leq\mc{H}(k_3+8,\rho',\up)$$
such that $G_2=\frac{H}{\Delta}\in\mc{H}(k_3-4,\rho',\up)$ has the form \eqr{arbF} with $n_{j_1}\geq1$. Now, setting
$(k,N)=(-2,3)$ in part two of Lemma \ref{lem:5dtdet} informs us that $\mc{H}(k_3-4,\rho',\up)$ is \emph{at most} two-dimensional,
and we claim that in fact $\mc{H}(k_3-4,\rho',\up)=\langle DG_1,G_2\rangle$ is \emph{exactly} two-dimensional. To see this,
we must show that no relation $DG_1=\alpha G_2$ exists, with $\alpha\in\bb{C}$. Writing $G_1$ in the form \eqr{arbF}
and examining the modular Wronskian $W(G_1)$, we find from the bound \eqr{wtbound} that $n_j=0$ for $1\leq j\leq5$ in \eqr{arbF}.
Furthermore, because of the way we chose $j_1$, it follows directly from the definition of the modular derivative that when
we write $DG_1$ in the form \eqr{arbF}, we again have $n_{j_1}=0$. On the other hand, by definition we know that $G_2$, when
written as $\eqr{arbF}$, satisfies $n_{j_1}\geq1$. Therefore $DG_1\neq\alpha G_2$ for any $\alpha\in\bb{C}$, and our claim
about $\mc{H}(k_3-4,\rho',\up)$ is verified. Part two of Lemma \ref{lem:5dtdet} with $(k,N)=(-1,3)$ and Lemma \ref{lem:dfind}
show that $\mc{H}(k_3-2,\rho',\up)=\langle D^2G_1,E_4G_1\rangle$ is two-dimensional, and this is enough information to
establish that $\mc{H}(\rho',\up)$ has the Hilbert-Poincar\'e series \eqr{hpfinal} with $N=3$, $n_3=-3$, as claimed.

$\textbf{N=4.}$ Set $F_4=F_{(1,1,1,1,0)}\in\mc{H}(k_4,\rho',\up)$ in \eqr{5dMMDEvec}. We find, as usual, a vector
$$G\in\langle D^2F_4,E_4F_4\rangle\leq\mc{H}(k_4+4,\rho',\up)$$
such that $G_1=\frac{G}{\Delta}$ is a nonzero vector in $\mc{H}(k_4-8,\rho',\up)$. Part one of Lemma \ref{lem:5dtdet} implies
that $k_4-8$ is the minimal weight for $\mc{H}(\rho',\up)$, and part two with $(k,N)=(-4,4)$ shows that
$\mc{H}(k_4-8,\rho',\up)=\langle G_1\rangle$ is one-dimensional. Similar to the $N=3$ case above, let us fix an $i_1\in\{1,\cdots,5\}$
such that $\lambda_{i_1}\neq\frac{k_4-8}{12}$. Using this $i_1$ to define the appropriate linear functionals, we may locate a vector
$$H\in\langle D^3F_4,E_4DF_4,E_6F_4\rangle\leq\mc{H}(k_4+6,\rho',\up)$$
such that $G_2=\frac{H}{\Delta}$ is a nonzero vector in $\mc{H}(k_4-6,\rho',\up)$ which, when written in the form \eqr{arbF},
has the property that $n_{i_1}\geq1$. Once again, it then follows directly from the definition of the modular derivative that
$DG_1$ and $G_2$ are linearly independent, and part two of Lemma \ref{lem:5dtdet} implies (using $(k,N)=(-3,4)$) that
$\mc{H}(k_3-6,\rho',\up)=\langle DG_1,G_2\rangle$ is two-dimensional. We iterate this logic one last time, defining
an $i_2\in\{1,\cdots,5\}-\{i_1\}$ such that $\lambda_{i_2}\neq\frac{k_4-6}{12}$, and using this $i_2$ we define the
appropriate linear functionals to produce a vector
$$\widetilde{G}\in\langle D^4F_4,E_4D^2F_4,E_6DF_4,E_8F_4\rangle\leq\mc{H}(k_4+8,\rho',\up)$$
such that $G_3=\frac{\widetilde{G}}{\Delta}$ is a nonzero vector in $\mc{H}(k_4-4,\rho',\up)$, with the property that, when
written in the form \eqr{arbF}, we have $n_j\geq1$ for $j=i_1,i_2$. Using $(k,N)=(-2,4)$ in part two of Lemma \ref{lem:5dtdet}
shows that $\mc{H}(k_4-4,\rho',\up)$ is at most three-dimensional, and we claim that this bound is realized, \ie
$$\mc{H}(k_4-4,\rho',\up)=\langle E_4G_1,DG_2,G_3\rangle.$$
To verify this, assume there is a relation
\begin{equation}\label{eq:fineq}
\alpha_1E_4G_1+\alpha_2DG_2+\alpha_3G_3=0
\end{equation}
for some $\alpha_j\in\bb{C}$. Now, it follows directly from the definition of $G_1$ and the bound obtained from Theorem
\ref{thm:MMDEwt} that when $G_1$ is written in the form \eqr{arbF}, we have $n_j=0$ for $1\leq j\leq5$. Similarly, when
$G_2$ is written in this form we have $n_{i_1}=1$ and $n_j=0$ otherwise, and for $G_3$ we find that $n_{i_1}=n_{i_2}=1$ and
$n_j=0$ otherwise. From this information, one sees that the $j^{th}$ component of the left-hand side of \eqr{fineq} is of
the form
\begin{eqnarray}
&\ &\left(\alpha_1a_j(0)q^{\lambda_j}+\cdots\right)+\nonumber\\
\nonumber\\
&\ &\left(\alpha_2b_j(0)\left(\lambda_j+\delta_{j,i_1}-\frac{k_4-6}{12}\right)
q^{\lambda_j+\delta_{j,i_1}}+\cdots\right)+\\
\nonumber\\
&\ &\left(\alpha_3c_j(0)q^{\lambda_j+n_j}+\cdots\right),\nonumber
\end{eqnarray}
where $a_j(0)q^{\lambda_j}$, $b_j(0)q^{\lambda_j+\delta_{j,i_1}}$, $c_j(0)q^{\lambda_j+n_j}$ denote the leading terms
of $G_1$, $DG_2$, $G_3$, respectively, and $n_j$ is as stated above for $G_3$. Setting $j=i_1$ in \eq{0}, we find
that $\alpha_1=0$, since $a_{i_1}(0)\neq0$ and the second and third expansions in \eq{0} have no $q^{\lambda_{i_1}}$ term.
Similarly, the $j=i_2$ version of \eq{0} shows that $\alpha_2=0$, since we now know that $\alpha_1=0$, the second expansion
in \eq{0} has leading term $q^{\lambda_{i_2}}$, and the third has leading term $q^{\lambda_{i_2}+1}$. But then $\alpha_3=0$
also, and this proves that $\mc{H}(k_4-4,\rho',\up)$ is three-dimensional. We now have have enough information
to see that $\{G_1,DG_1,G_2,DG_2,G_3\}$ forms a set of free generators for $\mc{H}(\rho',\up)$ as $\m$-module, so that the
Hilbert-Poincar\'e series for $\mc{H}(\rho',\up)$ is of the form \eqr{hpfinal}, with $N=4$, $n_4=-4$.
\end{proof}

Although we omit them here in favor of concision, we note that explicit $\m$-bases have been computed for the $N=1,2,3$ cases
of Theorem \ref{thm:finthm}, and may be found in \cite[Secs 4.5.2-4.5.4]{Marks1}. Similarly, one may find in \lc\ explicit
formulas for the dimensions of the various spaces $\mc{H}(k,\rho,\up)$; of course this information may be also be
obtained directly from the given Hilbert-Poincar\'e series and the classical formula for the dimension of $\m_k$.

\appendix
\section{Cusp forms and modular differential equations}

One way of seeing why Definition \ref{def:Tdet} is not frivolous in arbitrary dimension is to examine the effect of cusp
forms on the MMDE theory reviewed in Subsection \ref{sec:mlde}. For concreteness, we focus on dimension six. Because $\m_{12}=\bb{C}E_{12}\oplus\bb{C}\Delta$, one observes that the arbitrary MMDE of order six is of the form $(L+c\Delta)[f]=0$,
where $L$ is a sixth order Eisenstein operator \eqr{eop}, and $c$ an arbitrary complex number.
Now, Lemma \ref{lem:eop} states that for any set $\Lambda$ of six complex numbers, there is a unique Eisenstein operator
$L_\Lambda$, such that the MMDE $L_\Lambda[f]=0$ has the indicial roots $\Lambda$. Recall (\cf \cite{H}) that if
this MMDE is written in the form \eqr{rewrite}, the indicial polynomial (whose roots are the indicial roots of
 \eqr{rewrite}) is determined by the constant terms of the holomorphic functions $g_j(q)$. But for any $c\in\bb{C}$, $c\Delta$
has constant term 0. Thus we find that for each set $\Lambda$, every operator in the family
$\{L_\Lambda+c\Delta\ |\ c\in\bb{C}\}$ has $\Lambda$ as its set of indicial roots. Using this fact, it is easy to construct
irreducible representations of dimension six which are \emph{not} $T$-determined:

For example, let $\Lambda=\{r_1,\cdots,r_5\}$ be any set of distinct real numbers satisfying $0<r_j<1$ for each $j$,
such that no proper sub-sum of the $r_j$ is of the form $\frac{x}{12}$, $x\in\bb{Z}$, and such that $r=\sum r_j=\frac{5}{2}$.
(For example, one may choose $\Lambda=\{\frac{2}{22},\frac{5}{22},\frac{8}{22},\frac{19}{22},\frac{21}{22}\}$.) Using these $r_j$ as
indicial roots, one obtains by Corollary \ref{cor:ummde} a unique Eisenstein operator in weight (by Theorem \ref{thm:MMDEwt})
$\frac{12r}{5}-4=2$, of the form
$$L_\Lambda=D_2^5+\alpha_4E_4D_2^3+\cdots+\alpha_{10}E_{10}.$$
For each $c\in\bb{C}$, consider the operator
$$L_c=L_\Lambda D_0+c\Delta=D_0^6+\alpha_4E_4D_0^4\cdots+\alpha_{10}E_{10}D_0+c\Delta,$$
and the associated representation $\rho_c:\G\rightarrow\gln{6}$, arising from the $|_0$ action of $\G$ on the solution
space $V_c$ of the MMDE $L_c[f]=0$, as afforded us by Theorem \ref{thm:mldevec}. Recalling that the solution space of
$D_kf=0$ is spanned by $\eta^{2k}$, one sees immediately that any constant function will be a solution of $L_0[f]=0$, thus
the solution space $V_0$ contains a one-dimensional subspace $\m_0=\bb{C}$ of functions which are invariant under $|_0$.
In particular, the representation $\rho_0$ is indecomposable (by Lemma \ref{lem:indecMMDE}) but \emph{not} irreducible.

As mentioned above, the indicial roots of $L_c[f]=0$ will be the same for any $c\in\bb{C}$, and it is
clear that the indicial roots of $L_0[f]=0$ are $\{0,r_1,\cdots,r_5\}$: The solutions of $L_0[f]=0$ are exactly the functions
$\{f\ |\ L_\Lambda[D_0f]=0\}$, and if $f=q^\lambda+\cdots$, then $D_0f=q\frac{df}{dq}=\lambda q^\lambda+\cdots$, so either $\lambda=0$
or the solution $f$ satisfying $L_0[f]=0$ has the same leading exponent as the solution $g=D_0f$ satisfying $L_\Lambda[g]=0$.
In particular, for any $c\in\bb{C}$ we have, say, $\rho_c(T)=\diag{1,\e{r_1},\cdots,\e{r_5}}$.

We claim that $\rho_c$ is irreducible for any $c\in\bb{C}^\ast$. To see this, observe that for any $c$, the set $\{0,r_1,\cdots,r_5\}$
of indicial roots of $L_c$ has the same sub-sum property as $\Lambda$, with the obvious exceptions $0=\frac{0}{12}$ and
$\sum r_j=\frac{5}{2}$. But any proper invariant subspace of a solution space $V_c$ must correspond to some proper sub-sum of
indicial roots of the form $\frac{x}{12}$, since this subspace defines a sub-representation of $\G$. Thus for any $c\in\bb{C}$,
the only possibilities are that the proper invariant subspace of $V_c$ is one-dimensional, and corresponds to the single indicial root 0,
or that it is five dimensional, and corresponds to the sub-sum $\frac{5}{2}=\sum r_j$. In the former case, the invariant subspace
$V$ must again consist of holomorphic modular forms of weight 0 (since it must be invariant under $|_0$ and is one-dimensional), so that $V=\m_0=\bb{C}$ again consists of the constant functions. But clearly a constant function solves $L_c[f]=0$ exactly when $c=0$
(since these functions already satisfy $L_\Lambda D_0[f]=0$), thus the former case is impossible. In the latter case, the invariant
subspace $V$ yields a sub-representation $\rho:\G\rightarrow\gln{5}$ which is evidently irreducible, and has the property that $\rho(T)$ and
$\rho_\Lambda(T)$ have the same eigenvalues, where $\rho_\Lambda$ denotes any representation associated to the MMDE
$L_\Lambda[f]=0$. Thus $\rho\sim\rho_\Lambda$, by Theorem \ref{thm:tw}. But $\mc{H}(\rho_\Lambda,\textbf{1})$ is a cyclic
$\mc{R}$-module with minimal weight 2, by Theorem \ref{thm:cyclic}, whereas there is a nonzero vector
$F\in\mc{H}(0,\rho,\textbf{1})$, whose components form a basis of the invariant subspace $V$. This is a contradiction, in light of
the isomorphism \eqr{isom}. Thus the latter case is also ruled out, and our claim is verified, \ie $\rho_c$ is irreducible for each $c\in\bb{C}^\ast$.

Along the same line of reasoning, there is one more important aspect of this example which must be mentioned. By Theorem
\ref{thm:cyclic}, each space $\mc{H}(\rho_c,\textbf{1})$ is a cyclic $\mc{R}$-module $\mc{R}F_c$, where the minimal weight vector $F_c\in\mc{H}(0,\rho_c,\textbf{1})$ has components which span the solution space $V_c$ of the MMDE $L_c[f]=0$. In particular, each space $\mc{H}(0,\rho_c,\textbf{1})=\langle F_c\rangle$ is one-dimensional. Using this fact, it is easy to see that $\rho_{c_1}$ is equivalent to $\rho_{c_2}$ if, and only if, $c_1=c_2$. This follows directly from the isomorphism \eqr{isom}: If $\rho_{c_1}\sim\rho_{c_2}$, then there is a $U\in\gln{6}$ such that $UF_{c_1}=\alpha F_{c_2}$ for some $\alpha\in\bb{C}$, and this implies that every component of $F_{c_2}$ is a solution of
$L_{c_1}[f]=0$ and vice-versa. Clearly this happens if, and only if, $c_1=c_2$.

Finally, we point out that this example obviously generalizes to any dimension/order greater than 5, since there will always be
cusp forms available to be utilized in this same manner. Indeed, it is hoped that a deep understanding of this type of example may
eventually lead to some sort of progress in solving the general problem of classifying irreducible (or, even better, indecomposable)
representations of $\G$ in arbitrary dimension.

We summarize the above discussion by recording the
\begin{prop}
In each dimension $d\geq6$, there exists a one parameter family of inequivalent indecomposable, $T$-unitarizable
representations
$$\{\rho_c:\G\rightarrow\gln{d}\ |\ c\in\bb{C}\}$$
with the following properties:
\begin{enumerate}
\item For every $c_1,c_2\in\bb{C}$, $\rho_{c_1}(T)=\rho_{c_2}(T)$.\\
\item $\rho_c$ is irreducible if, and only if, $c\in\bb{C}^\ast$.
\end{enumerate}
\qed
\end{prop}

This Proposition gives an indication of just how spectacularly the results of \cite{TW} fail to be true in dimension greater
than five; see also Remark 2.11.3, \lc.

\subsection*{Acknowledgement} This paper is a mild strengthening of the main results of the author's doctoral dissertation, written
under the supervision of Geoffrey Mason at the University of California, Santa Cruz. We are very happy to acknowledge the University
-- and in particular Professor Mason -- for the support and encouragement we received during our time there. We would also like to thank the
Max-Planck-Institut f\"ur Mathematik in Bonn for providing a truly glorious work environment during our visit there, which consequently
made the preparation of this document a much more enjoyable task. Thanks also to Terry Gannon and Martin Weissman, who read and advised
on earlier versions of this paper.

\bibliographystyle{amsplain}
\bibliography{modularbib}
\end{document}